\documentclass{article}
\usepackage{amsmath}
\usepackage{amssymb}

\newcommand{\R}{\mathbb{R}}

\newcommand{\Z}{\mathbb{Z}}
\newcommand{\Zp}{\mathbb{Z}_{p}}

\newcommand{\Qp}{\mathbb{Q}_{p}}
\newcommand{\GL}{\mathrm{GL}}
\newcommand{\Gm}{\mathbb{G}_{\mathrm{m}}}

\newtheorem{prop}{Proposition}
\newtheorem{theo}{Th\'eor\`eme}
\newtheorem{cor}{Corollaire}
\newtheorem{lemme}{Lemme}

\title{Existence de $F$-cristaux avec structures suppl\'ementaires.}

\author{J.-P. ~\textsc{Wintenberger}}

\begin{document}

\maketitle

Dans cet article, nous prouvons, sous des hypoth\`eses
de non ramification,  une conjecture de Kottwitz
et Rapoport sur l'existence de cristaux avec structures
suppl\'ementaires (\cite{[KR02]}, \cite{[R02]}).

Essayons d'expliquer ce dont il s'agit.

Soient $p$ un nombre premier, $k$  un corps alg\'ebriquement 
clos de caract\'eristique $p$, 
$W(k)$  l'anneau des vecteurs de Witt \`a coefficients
dans $k$. Dans cette introduction, soient
$L=\Qp \otimes_{\Zp} W(k)$, et $\overline{L}$
une cl\^oture alg\'ebrique de $L$.
Notons $\sigma$ le Frobenius sur $k$, $W(k)$, et $L$. 
Un ($F$-)cristal $M$  est ici un $W(k)$-module libre de 
type fini et une application  $ \sigma$-lin\'eaire injective $\phi$
de $M$ dans $L \otimes_{W(k)} M$. 
Dieudonn\'e a 
a classifi\'e les cristaux \`a isog\'enie pr\`es \it i.e. \rm les 
$L$-espaces vectoriels de dimension finie munis d'une application 
$\sigma$-lin\'eaire bijective.
Les isocristaux sont 
classif\'es par leur polygone de Newton.
Mazur a introduit pour les cristaux 
un autre invariant, le polygone de Hodge, qui 
est donn\'e par les diviseurs \'el\'ementaires de $\phi (M)$ 
par rapport \`a $M$ (\cite{[M72]}). 
Mazur a observ\'e que le polygone de Hodge est en 
dessous du polygone de Newton, et que les deux polygones  ont m\^emes 
extr\'emit\'es. Kottwitz et Rapoport ont pos\'e (et r\'esolu) le probl\`eme
d'une r\'eciproque \`a l'observation de Mazur, \it i.e. \rm
l'existence d'un $F$-cristal de polygones de Hodge et Newton 
donn\'es, ces polygones v\'erifiant la condition de Mazur (\cite{[KR02]}).

Soit $G$ un groupe r\'eductif (connexe) sur $\Qp$. Kottwitz a d\'efini
la notion de $G$-isocristal  (\cite{[K85]},
\cite{[K97]}). 
Un tel isocristal est donn\'e par un \'el\'ement $b\in G( L)$ ;
on a $\phi=b\sigma$.
Deux \'el\'ements $b$ et $b'\in G(L)$ d\'efinissent  des $G$-isocristaux 
isomorphes si $b$ et $b'$ sont $\sigma$-conjugu\'es
par un \'el\'ement de $G( L)$.  La classe
de $\sigma$-conjugaison d'un \'el\'ement $b$ de $G(L)$ est  not\'ee $[b]$ et
l'ensemble des classes de $\sigma$-conjugaison  $B(G)$. 
Dans ce contexte, un cristal
devient la donn\'ee d'un parahorique $\widetilde{K}$ de $G(L)$. 
Le polygone de Hodge devient la position relative de 
$\phi (\widetilde{K})$ par rapport \`a $\widetilde{K}$.
Cette position relative est mesur\'ee par un invariant dans 
un ensemble de doubles classes 
$\widetilde{W}^K\backslash \widetilde{W}/\widetilde{W}^K$.
Le groupe $\widetilde{W}$ est le produit semi-direct
du groupe de Weyl affine de l'immeuble de $G_L$ par le groupe 
ab\'elien de type fini $\pi_1 (G)_{\mathrm{Gal}(\overline{L}/L)}$
et $\widetilde{W}^K$ le groupe de Weyl de $\widetilde{K}$
(\cite{[R02]},\cite{[B98]}). 
Pour $K$ un parahorique de $G$, $\overline{w}\in \widetilde{W}^K\backslash 
\widetilde{W}/\widetilde{W}^K$ et $b\in B(G)$, 
Kottwitz et Rapoport d\'efinissent un sous-ensemble
$X_{\overline{w}}(b)_{K}$ qu'ils appellent vari\'et\'e de 
Deligne-Lusztig affine. 
Grossi\`erement, $X_{\overline{w}}(b)_{K}$ est l'ensemble
des parahorique $\widetilde{K'}$ de $G_L$ qui sont conjugu\'es de $K$
et tels que la position de $\phi (\widetilde{K}' )$ par rapport
\`a $\widetilde{K}'$ soit $\overline{w}$. 
Le probl\`eme de la 
r\'eciproque au th\'eor\`eme de Mazur devient celui
de d\'eterminer pour quels $[b]$ et 
$\overline{w}$,  $X_{\overline{w}}(b)_{K}$ est non vide.

Supposons $G$ non ramifi\'e sur $\Qp$.
Si $K$ est hypersp\'ecial, et $\widetilde{K}$ le parahorique
de $G_L$ obtenu par extension des scalaires de $\Qp$ \`a $L$,
$\widetilde{W}^K\backslash \widetilde{W}/\widetilde{W}^K$
s'identifie \`a l'ensemble des classes de conjugaison de 
groupes \`a un param\`etre de $G_L$. 
Pour $\{ \mu \}$ une telle classe de conjugaison, 
Kottwitz a d\'efini un sous-ensemble $B(G, \{ \mu \})$ de 
$B(G)$ (\cite{[K97]}). La condition d\'efinissant $B(G, \{ \mu \})$
g\'en\'eralise celle de Mazur. Rapoport et Richartz ont prouv\'e 
que si $X_{\overline{w}}(b)_{K}$ est non vide, $[b]\in
B(G, \{ \mu \})$ (\cite{[RR96]}). Kottwitz et Rapoport
conjecture la r\'eciproque (conj. 4.6. de \cite{[R02]} ou 
\cite{[KR02]}).
Ils la prouvent en particulier pour
les groupes lin\'eaires $\GL_{n}$ et symplectiques (\cite{[KR02]}).
Leigh la prouve pour les groupes classiques d\'eploy\'es (\cite{[Le03]}).
Nous la prouvons pour $\{ \mu \}$ minuscule (cor. \ref{conjhyper}).

Soient  $K$ un parahorique de $G$ non n\'ecessairement hypersp\'ecial
et $[b]\in B(G)$. Il semble difficile de d\'ecrire, m\^eme 
conjecturalement, quels sont les 
$\overline{w}\in \widetilde{W}^K\backslash \widetilde{W}/\widetilde{W}^K$  
qui sont tels que $X_{\overline{w}}(b)_{K}$
soit  non vide, bien que ceci ait \'et\'e fait dans des 
cas particuliers (exemple 4.3. de      \cite{[R02]}, \cite{[Re03]}). 
Soit $\{ \mu \}$ une classe de conjugaison de groupes \`a
un param\`etre de $G_{\overline{L}}$. Kottwitz et Rapoport 
d\'efinissent un sous-ensemble $\mathrm{Adm}(\{ \mu \} )_{K}$ de 
$\widetilde{W}^K\backslash \widetilde{W}/\widetilde{W}^K$. 
La r\'eunion des  vari\'et\'es de 
Deligne-Lusztig affines  $X_{\overline{w}}(b)_{K}$ pour 
$\overline{w}\in \mathrm{Adm}(\{ \mu \} )_{K}$ est not\'ee
$X(\{ \mu \},b)_{K}$. Ils conjecturent :

- 1) $X(\{ \mu \},b)_{K}$ est non vide si 
et seulement si $[b]\in B(G, \{ \mu \})$ ;

- 2) si $K\subset K'$ l'application naturelle
$X(\{ \mu \},b)_{K}\rightarrow X(\{ \mu \},b)_{K'}$
est surjective ;

(en fait ils demandent  si l'hypoth\`ese $\{ \mu \}$ minuscule
n'est pas n\'ecessaire pour la conjecture)(conj. 5.2. de \cite{[R02]} 
ou \cite{[KR02]}). 
Nous prouvons la partie 1) de la conjecture si $G$ est non ramifi\'e 
sans l'hypoth\`ese $\{ \mu \}$ minuscule
(th. \ref{conjadm}). 
Les parties \ref{posiwahori} et \ref{bonneposition} disent 
quelque chose sur la partie 2) de la conjecture ; cependant, faute
d'avoir d\'egag\'e le bon cadre, nous ne donnons pas d'\'enonc\'e.
Puisque, si $K$ est hypersp\'ecial et $\{ \mu \}$ minuscule,
l'ensemble  $\mathrm{Adm}(\{ \mu \} )_{K}$ est r\'eduit \`a 
$\{ \mu \}$ (cor. 3.12. de \cite{[R02]}),   la conjecture  
4.6. de \cite{[R02]} dans le cas $\{ \mu \}$ minuscule
est prouv\'ee comme une cons\'equence  de la partie 1)
de la conjecture 5.2..

Donnons quelques id\'ees des d\'emonstrations. 

Dans le cas o\`u $\{ \mu \}$ est minuscule et
$G$ est $\GL _n$ ou le groupe des similitudes symplectiques,
on peut prouver la conjecture 4.6. en prenant pour 
matrices de $\phi$ dans une base des  r\'eseaux cherch\'es
des matrices monomiales, mais le cas de $\GL_2$
 montre d\'ej\`a que cela n'est pas possible en g\'en\'eral si
$\{ \mu \}$ n'est pas minuscule. 
En fait, m\^eme dans les cas de $\GL _{n}$, 
la preuve de Kottwitz et Rapoport est
beaucoup plus \'evolu\'ee : 
elle repose sur  une propri\'et\'e de positivit\'e de l'isomorphisme
de Satake (\cite{[R01]}). 

Fontaine et Rapoport ont donn\'e une autre
preuve du th\'eor\`eme de Kottwitz et Rapoport (\cite{[FR02]}).
Ils prouvent, sous l'une des deux conditions
qui d\'efinissent l'ensemble $B(G,\{ \mu \})$, 
l'existence d'une filtration admissible d\'efinie par un groupe \`a
un param\`etre $\mu \in \{ \mu \}$. L'admissibilit\'e est \`a 
comprendre dans le sens la th\'eorie
de Fontaine  des modules de Dieudonn\'e filtr\'es admissibles
(\cite{[FL83]},\cite{[F94]},\cite{[FC01]},\cite{[Fa95]}) .
Un th\'eor\`eme
de Laffaille donne alors l'existence d'un r\'eseau $M$
de $D$ qui est fortement divisible (\cite{[L80]}).
Dans le cas de $\GL _n$, un tel r\'eseau convient 
pour la conjecture 4.6.. 
 
Sous de plus la deuxi\`eme
condition d\'efinissant $B(G,\{ \mu \})$, nous prouvons l'existence 
d'un parahorique $K$ de $G$ (d\'efini sur $F$) tel que la conjecture 
5.2. de Kottwitz et Rapoport soit vraie pour $K$ (th. \ref{bonpara}).
Donnons une raison heuristique pour laquelle ceci peut
\^etre vrai. Si
$\mu$ d\'efinit une filtration admissible, pour toute repr\'esentation
lin\'eaire de $G$, le th\'eor\`eme de Colmez et Fontaine permet
d'associer une repr\'esentation $p$-adique cristalline de 
$\mathrm{Gal}(\overline{L}/L )$. Les r\'eseaux fortement 
divisibles correspondent aux r\'eseaux stables par 
$\mathrm{Gal}(\overline{L}/L )$. 
Comme $\mathrm{Gal}(\overline{L}/L )$ est compact,
son image dans $G(\Qp )$ est contenue dans un  parahorique
auquel correspond le parahorique $K$, par le foncteur de Fontaine entre
repr\'esentations $p$-adiques et modules filtr\'e. 
Nous prouvons dans \ref{pargal} la version galoisienne
du th\'eor\`eme \ref{bonpara}. 
Dans la partie \ref{parmdf}, nous prouvons le th\'eor\`eme 
\ref{bonpara} en travaillant du c\^ot\'e module de Dieudonn\'e
filtr\'e et en utilisant que la filtration dans 
la cat\'egorie de ces modules peut \^etre scind\'ee
de mani\`ere fonctorielle (\cite{[W84]}).
Dans la proposition 
\ref{calculinv}, nous calculons l'invariant 
dans $\widetilde{W}$ pour des Iwahori qui sont associ\'es
\`a certaines alc\^oves dont l'adh\'erence contient la facette
associ\'ee \`a $K$. Nous introduisons la notion de bonne position
pour ces alc\^oves. 
Nous prouvons qu'il existe des alc\^oves en bonne position
(prop. \ref{existebonneposition}).
Elle entra\^{\i}ne la conjecture 5.2.
pour les Iwahori associ\'es \`a ces alc\^oves 
(prop. \ref{bonnepositionbruhat}). Comme tous les 
Iwahori sont conjugu\'es, la conjecture en r\'esulte.  

On notera que la m\'ethode par laquelle   nous raffinons 
la $K$-structure en une structure iwahorique est 
tr\`es li\'ee aux structures alg\'ebriques apparaissant
dans \cite{[GP68]}, \cite{[Kraft]} et \cite{[Mo01]}
(remarque \ref{GP}).
D'autre part, nous prouvons en fait un peu
plus que la conjecture 5.2. : nous prouvons 
son analogue lorsqu'on remplace dans la d\'efinition
de l'ensemble $\mathrm{Adm}(\{ \mu \} )$ l'ordre de 
Bruhat par l'ordre de Bruhat faible (la terminologie
est que si $w'$ est  inf\'erieur pour l'ordre de Bruhat
faible \`a $w$, $w'$ l'est pour l'ordre de Bruhat usuel)
(remarque \ref{Bruhatfaible}).  

Je remercie Rapoport pour son amical int\'er\^et.
Je remercie Verma pour m'avoir dit qu'il devait \^etre
possible de comparer $t_{\lambda}$ et 
$t_{\lambda}\theta$ (prop. \ref{bonnepositionbruhat}),
Sujatha et Parimala pour m'avoir aid\'e \`a retrouver
la r\'ef\'erence \cite{[BT71]}, et Haines pour ses 
remarques sur le \ref{G'}. L'utilisation que nous faisons de 
\cite{[BT71]} est inspir\'ee du 3 de \cite{[G98]}.

Nous travaillons en fait dans un cadre un peu plus large.
Soit  $F$ est une 
extension finie de $\Qp$ (ci-dessus $F=\Qp$). On suppose que le corps
alg\'ebriquement clos $k$  
contient  le corps r\'esiduel $k_F$ de $F$.
On pose :
$L=F\otimes_{W(k_F)} W(k)$. 
On d\'esigne toujours par $\overline{L}$ une cl\^oture alg\'ebrique
de $L$. On note $q$ le cardinal de $k_F$, $\sigma_0$
l'automorphisme de $W(k)$ induit par fonctorialit\'e par l'automorphisme
$x\mapsto x^q$ de $W(k)$ et $\sigma$ l'automorphisme 
$\mathrm{id}\otimes \sigma_0$ de $L$.
On d\'esigne par $O_F$ (resp. $O_L$) l'anneau de 
valuation de $F$ (resp. $L$).
On d\'esigne par $G$
un groupe r\'eductif (connexe) sur $F$. On d\'esigne par
$G_{\mathrm{ad}}$ le groupe adjoint de $G$, et par
$\mathcal{B}(G_{\mathrm{ad}},.)$ son immeuble sur 
$F$ ou $L$ selon que le point est $F$ ou $L$.

\section{Existence de structures parahoriques :
le c\^ot\'e galoisien.}\label{pargal}

\subsection{}\label{11} Soit $I_L$ le groupe de Galois de $\overline{L}/L$. 
Soit 
$\rho : I_L\rightarrow G(F)$ une repr\'esentation 
semi-stable, autrement dit, pour toute
repr\'esentation lin\'eaire $\rho _{U}$ de $G$
dans un $F$-espace vectoriel de dimension finie $U$,
la repr\'esentation $p$-adique 
$\rho _{U}\circ \rho$ semi-stable
en tant que repr\'esentation $p$-adique 
dans le $\Qp$-espace vectoriel $U$
(au sens de \cite{[F94]}). On dit que $\rho$
v\'erifie $(*)$  
si pour toute repr\'esentation $\rho_U$, la repr\'esentation
de $I_L$ dans $U$ obtenue par composition de $\rho$ et 
de $\rho_U$ v\'erifie la condition $(*)$ de 
la remarque 3 de \cite{[FR02]} (voir aussi 1.40 de \cite{[RZ96]}).
On v\'erifie facilement qu'il  suffit pour ceci que ce soit vrai 
pour une repr\'esentation fid\`ele $\rho_U$. La condition
est vide si $F=\Qp$. Le foncteur de Fontaine \'etablit
une \'equivalence entre ces repr\'esentations galoisiennes
et les $(\phi ,N)$-modules  filtr\'es admissibles sur $L$
(\cite{[FC01]}, la d\'efinition d'un
$\phi$ module filtr\'e sur $L$ est rappel\'ee au \ref{paireadm}).

\begin{prop}\label{pargalprop} Soit  
$\rho : I_L\rightarrow G(F)$ une repr\'esentation
$p$-adique. On suppose qu'elle est semi-stable et 
qu'elle v\'erifie la condition $(*)$.
Alors, il existe un parahorique $K$ de $G$ qui contient 
l'image $\rho (I_L )$ du groupe de Galois.
\end{prop}

\emph{D\'emonstration.} Le groupe $\rho (I_K )$ est compact.
Il  fixe donc un point $x$ dans l'immeuble $\mathcal{B}(G_{\mathrm{ad}},F)$
(2.3.1. de \cite{[T79]}). Les points de l'immeuble 
$\mathcal{B}(G_{\mathrm{ad}},F)$ sont les points fixes
par $\sigma$ de l'immeuble $\mathcal{B}(G_{\mathrm{ad}},L)$ 
(5.1.25 et 4.2.25 de \cite{[BT84]}).
Le point $x$ 
appartient \`a une unique  facette
$\mathcal{F}$ de $\mathcal{B}(G_{\mathrm{ad}},L)$. Cette facette est 
stable par $\sigma$.
Elle 
d\'efinit un parahorique $K_{\mathcal{F}}\subset G(F)$
de $G$ (5.2.6. de \cite{[BT84]}).
D'apr\`es \cite{[HR?]},
ce parahorique $K_{\mathcal{F}}$
peut \^etre d\'efini comme l'intersection du stabilisateur dans $G(F)$
de la facette $\mathcal{F}$ avec le noyau de l'homomorphisme
$\widetilde{\kappa}_G : G(L)\rightarrow \pi_1 (G)_{I_L}$ d\'efini
par R. Kottwitz dans le paragraphe 7 de \cite{[K97]}
(voir aussi \cite{[B98]} pour la d\'efinition
de $\pi _1 (G)$).
Ceci r\'esulte de ce que le noyau de $\widetilde{\kappa}_G$
est le sous-groupe $G(L)'$ de $G(L)$ d\'efini au 5.2.11 de \cite{[BT84]}
(voir le \ref{G'} ).
La proposition r\'esulte alors du  lemme :

\begin{lemme}\label{kgal} Soit $\rho : I_L \rightarrow G(F)$
une $F$-repr\'esentation qui est semi-stable. Alors, $\rho(I_L )$
est contenue dans le noyau de l'homomorphisme de Kottwitz
$\widetilde{\kappa}_G$.
\end{lemme}

\emph{D\'emonstration.} On sait (3.1.2 de \cite{[W97]}
et introduction de \cite{[FC01]}) qu'il existe :

- un groupe r\'eductif $G'$ sur $F$ dont le sous-groupe des 
commutateurs $[G',G']$ est simplement connexe ;

- un morphisme $f :G'\rightarrow G$ dont le noyau est ab\'elien 
et une $F$-repr\'esentation semi-stable 
$\rho' : I_L\rightarrow G' (F)$ telle que $f\circ \rho'=\rho$.

A dire vrai, ceci n'est \'enonc\'e dans les articles
cit\'es que lorsque $F=\Qp$, mais
la d\'emonstration qui s'y trouve marche pour $F$ quelconque.

On peut de plus supposer que $G'/[G',G']$ est un tore induit.
En effet, soit pour $F'$  extension finie
de $F$ contenue dans $\overline{L}$,
$T_{F'}$ le tore obtenu par restriction des scalaires 
\`a la Weil de $F'$ \`a $F$ \`a partir du groupe
multiplicatif sur $F'$, et soit $\rho_{F'} : I_{L}\rightarrow 
T_{F'}(F)={F'}^*$ la repr\'esentation $p$-adique associ\'ee
\`a un groupe formel de Lubin-Tate pour $F'$ (2 de
\cite{[S79]}).
La repr\'esentation
$p$-adique $\rho_{F'}$ est une $F$-repr\'esentation et elle
est cristalline. Soit $F'$ telle que $G'/[G',G']$ soit d\'eploy\'e
sur $F'$. La repr\'esentation $I_L\rightarrow G'(F )\rightarrow
(G'/[G',G'])(F)$ est un objet de la sous-cat\'egorie tannakienne
de la cat\'egorie des $F$-repr\'esentation engendr\'ee par $\rho_{F'}$
(Serre \it loc. cit.\rm).
Ceci d\'efinit un morphisme $f''$ de $T_{F'}$ dans $G'/[G',G']$.
On remplace $G'$ par le produit fibr\'e de $G'$ et $T_{F'}$
au dessus de $G'/[G',G']$, $G'\rightarrow G'/[G',G']$ 
\'etant la projection et $T_{F'}\rightarrow G'/[G',G']$ 
\'etant $f''$. 

Comme le morphisme $\widetilde{\kappa}_{G}$ est fonctoriel
en $G$ (7.4 de \cite{[K97]}), il suffit de prouver la 
proposition pour $\rho '$. Comme $[G',G']$ est simplement 
connexe, 
il suffit de la prouver pour le compos\'e de 
$\rho '$ avec la projection  $G'(F)\rightarrow (G'/[G',G'])(F)$.
Comme $T:=G'/[G',G']$ est un tore induit, le noyau de 
$\widetilde{\kappa}_T$ est le sous-groupe born\'e maximal
de $T(L)$, et alors la proposition r\'esulte de la compacit\'e
de $I_L$.

\subsection{Remarque.}\label{12}M\^eme si $G$ est suppos\'e
non ramifi\'e, l'exemple suivant montre que l'on ne peut pas 
toujours supposer le parahorique $K$ hypersp\'ecial. 
Cet exemple contredit la proposition 5.31 de \cite{[RZ96]}. 

On suppose $F=\Qp$.
Soient $E$ l'extension quadratique non ramifi\'ee de $\Qp$
contenue dans $\overline{L}$,
$O_E$ l'anneau de ses entiers et $U_E$ le groupe des 
unit\'es de $O_E$. Soient $\Lambda$ un $O_E$-module libre 
de rang $2$, $\Lambda =\Lambda_1 \oplus \Lambda_2$ avec
$\Lambda_1$ et $\Lambda_2$ libres de rang $1$ sur $O_E$.
Soit $\rho_E : I_L \rightarrow U_E $ la repr\'esentation
galoisienne associ\'ee \`a un groupe formel de Lubin-Tate
pour $E$.
On fait agir $I_L$ sur $\Lambda$ par $\rho_E$
suivi de l'action de $U_E$ par homoth\'eties. 
On a donc la repr\'esentation galoisienne $\rho$ associ\'ee
\`a un groupe formel produit de deux groupes formels de Lubin-Tate
; $U:=\Qp \otimes_{\Zp} \Lambda$. 
On choisit des bases de $\Lambda_1$ et $\Lambda_2$,
en tant que $\Zp$-modules, et donc des formes altern\'ees
$\mathrm{det}_1$ et $\mathrm{det_2}$ pour lesquelles
$\Lambda_1$ et $\Lambda_2$ sont auto-duaux. 
On d\'efinit la  forme symplectique $(\ ,\ )$ par
$(\ ,\ )=p \mathrm{det}_1 \oplus \mathrm{det}_2$.
Le  groupe $G$ est le groupe des similitudes
symplectiques. La repr\'esentation galoisienne
$\rho$ est \`a valeurs dans $G(\Qp )$, 
le facteur de similitude \'etant donn\'e par le 
caract\`ere cyclotomique.
Le r\'eseau dual $\Lambda^*$ de $\Lambda$ est 
$p^{-1}\Lambda_1 \oplus \Lambda _2$, qui est aussi
stable par Galois. Tout r\'eseau $\Lambda '$ de la forme
$\Lambda'_1 \oplus \Lambda _2$, avec 
$\Lambda_1\subsetneqq \Lambda '_1 \subsetneqq \Lambda_1 ^*$
est autodual. Aucun n'est stable par Galois
puisque $\Lambda^* / \Lambda$ est isomorphe 
en tant que module galoisien \`a $O_E / pO_E$
avec action de Galois donn\'ee par $\rho _E
\mathrm{mod}\ p$ (ceci contredit l'assertion 
p.252 de \cite{[RZ96]} ; dans la derni\`ere ligne
de cette page, $\mathcal{W}$ est seulement 
isomorphe \`a un twist de $\mathcal{W}^*$). 

Il n'existe pas de $\Zp$ r\'eseau  $\Lambda ''$
de $U$ stable par $I_L$  qui soit homoth\'etique 
\`a son dual.  En effet,
pour un $\Zp$-r\'eseau $\Lambda ''$ de $U$ ,
on d\'efinit $o(\Lambda '' )$ comme \'etant la valuation du d\'eterminant
d'une base de $\Lambda ''$, \'evalu\'ee par rapport \`a une 
base d'un  r\'eseau autodual $\Lambda '$ comme ci-dessus.
Donc $o(\Lambda )= 1$, $o(\Lambda^*  )=-1 $.
Pour tout $\Zp$-r\'eseau $\Lambda ''$, on a $o({\Lambda''}^* )=
-o(\Lambda'' )$, donc si ${\Lambda''}^*$ est homoth\'etique
\`a $\Lambda''$, on a : $2o(\Lambda '')\equiv  0\ \mathrm{mod.}4$
et $o(\Lambda '')$ est paire. 
Si $\Lambda ''$ est stable par $I_L$, c'est un sous-$O_E$-module ;
comme $\mathrm{det}_{\Zp } = N_{E/ \Qp}\circ \mathrm{det}_E$,
$o(\Lambda '')$ a la m\^eme parit\'e que   
$o(\Lambda )$, donc est impaire.

\section{Existence de structures parahoriques :
le c\^ot\'e module de Dieudonn\'e filtr\'e.}\label{parmdf}

\subsection{Paires $(b,\mu)$ admissibles (\cite{[FR02]}).}
\label{paireadm}

Rappelons qu'un module de Dieudonn\'e filtr\'e $D$ (sur $L$)
est la donn\'ee d'un $L$-espace vectoriel de dimension
finie $D$, muni d'une application $\sigma$-lin\'eaire et bijective
$\phi$ et
d'une filtration $\mathcal{F}^{\bullet}_D$
par des sous-espaces vectoriels, filtration d\'ecroissante,
exhaustive et s\'epar\'ee.

Soit $b\in G(L)$ et $\mu : \Gm\rightarrow G$ un groupe \`a un 
param\`etre d\'efini sur $L$. Une repr\'esentation lin\'eaire
$\rho_U : G\rightarrow \GL _U$ de $G$ dans un $F$-espace vectoriel
$U$ de dimension finie d\'efinit 
un   module de Dieudonn\'e filtr\'e $D(b,\mu , \rho)$ dont l'espace 
vectoriel sous-jacent est $L\otimes_F U$ :

- le Frobenius $\phi$  est $ \rho_U (b)\circ (\sigma\otimes \mathrm{id}_U)$ ;

- la filtration est la filtration d\'ecroissante d\'efinie
par le groupe \`a un param\`etre $\rho_U\circ \mu$ :
si $D_{\mu}^{\bullet}$ est la graduation de $D=D(b,\mu , \rho)$
d\'efinie par $\mu$,
$\mathcal{F}^i _D =\sum_{i'\geq i} D_{\mu}^{i'}$. 

Rappelons que $(b,\mu)$ est \emph{admissible} si pour toute repr\'esentation
lin\'eaire $\rho$ le module de Dieudonn\'e filtr\'e $D(b,\mu ,\rho )$ 
est admissible (\cite{[FR02]}). 

Soit $P_{\mu}$ le parabolique de $G_L$ qui fixe la filtration
ci-dessus \it i.e. \rm $g\in G$ appartient \`a $P_{\mu}$ si, pour
toute repr\'esentation lin\'eaire $\rho$, $\rho (g)$ respecte la 
filtration de $D(\rho)$. Si $\mu '$ est un autre 
groupe \`a un param\`etre, on dit que $\mu$ et $\mu '$ \emph{d\'efinissent 
la m\^eme filtration}, si pour tout $\rho_U$ comme ci-dessus, les filtrations 
d\'ecroissantes sur $L\otimes_F U$  d\'efinies par  
$\rho_U\circ \mu$ et $\rho_U\circ \mu '$ co\"{\i}ncident, autrement 
dit s'il existe $u$ dans le radical unipotent
$ \mathrm{rad_u} (P_{\mu })(L)$ de $P_{\mu}$ tel que 
$\mu ' = \mathrm{int}(u)(\mu )$. Il est clair que 
si $\mu '$ d\'efinit la m\^eme filtration que $\mu$,
pour toute rep\'esentation lin\'eaire $\rho$ de $G$,
les modules de Dieudonn\'e filtr\'es $D(b,\mu , \rho)$
et $D(b,\mu ', \rho)$ co\"{\i}ncident.

Soit $(b_1,\mu_1)$ comme $(b,\mu )$, donc $b_1 \in G(L)$
et $\mu_1$ est un groupe \`a un param\`etre de $G_L$.
Soit $g\in G(L)$.
On dit que $(b_1,\mu_1)$ est \emph{le conjugu\'e par $g$ } de $(b,\mu )$
si $b_1=gb\sigma(g^{-1})$
et $\mu_1=\mathrm{int}(g)(\mu )$. 

\begin{prop}\label{conj} Soient $(b,\mu )$ et $(b_1,\mu_1)$ 
conjugu\'es par $g$.
Alors, pour 
toute repr\'esentation lin\'eaire $\rho$ de $G$ dans 
un $F$-espace vectoriel de dimension finie $U$, $\rho(g)$
induit un isomorphisme  des modules de Dieudonn\'e
filtr\'es  $D(b,\mu , \rho)$ et  
$D(b_1 ,\mu_1 , \rho)$. Si $(b,\mu )$
est admissible, il en est de m\^eme de $(b_1,\mu_1)$. \end{prop}

\emph{D\'emonstration.} Clairement,  la seconde assertion r\'esulte 
de la premi\`ere. 
Si $\phi$ et $\phi _1$ sont les Frobenius de 
$D(b,\mu , \rho)$ et $D(b_1 ,\mu_1 , \rho)$ respectivement,
on a :

$$\phi _1 \rho (g) = \rho(gb\sigma (g^{-1}))(\sigma\otimes 
\mathrm{id})\rho (g).$$ 

Comme : 

$$\rho(\sigma (g^{-1}))=(\sigma\otimes\mathrm{id})\rho (g^{-1} )
(\sigma^{-1}\otimes\mathrm{id}), $$

on voit que :

$$\phi _1 \rho (g)= \rho (g )\rho (b)(\sigma\otimes \mathrm{id})=
\rho(g )\phi.$$

De plus, comme $\mu_1=\mathrm{int}(g)(\mu )$, la filtration
de $D(b_1 ,\mu_1 ,\rho )$ est l'image de celle de 
$D(b ,\mu ,\rho )$ par $\rho (g)$, ce qui prouve la 
proposition.

\subsection{Paires $(b,\mu )$ super-admissibles.}\label{propsuper} 
Rappelons que si $D$ est un module de Dieudonn\'e filtr\'e sur $L$,
un $O_L$-r\'eseau $M$ de $D$ est dit fortement divisible, 
si l'on a : 

$$M= \sum_{i\in \Z} \pi^{-i}\phi ( \mathcal{F}^{i} \cap M),$$

$\pi$ \'etant une uniformisante quelconque de $F$.
Fontaine et Laffaille ont d\'efini une cat\'egorie 
de modules de Dieudonn\'e filtr\'es sur $O_F$ que nous notons
$\mathrm{MF}_{O_F}^{\mathrm{a}}$. Si on choisit $\pi$,
un objet est un $O_F$-module de type fini $M$
muni d'une filtration d\'ecroissante exhaustive et 
s\'epar\'ee $\mathcal{F}^{i}_M$ par des sous-modules
qui sont facteurs directs et d'applications 
$\sigma$-lin\'eaires $\phi^i : \mathcal{F}^{i}_M\rightarrow M$
telles que :

$$\phi^i_{\mid \mathcal{F}^{i+1}}= \pi \phi^{i+1},
\ M=\sum_{i\in \Z} \phi^i (\mathcal{F}^{i}_M ).$$ 

Un r\'eseau fortement divisible $M$ d'un module de Dieudonn\'e 
filtr\'e sur $L$ d\'efinit un objet de $\mathrm{MF}_{O_F}^{\mathrm{a}}$ :
on pose $\mathcal{F}^{i}_M=\mathcal{F}^{i}_D\cap M$ et 
$\phi^i$ est la restriction de  $ \pi^{-i}\phi$
\`a $\mathcal{F}^{i}_M$. 

Nous disons qu'un couple admissible $(b,\mu)$ est 
\it super-admissible \rm si l'on a :

- 1) $b=\sigma (\mu(p))$ : ceci \'equivaut \`a ce  que pour toute
repr\'esentation lin\'eaire $\rho$ de $G$, le Frobenius 
de $D(b,\mu ,\rho)$ soit $(\sigma\otimes \mathrm{id})\circ \mu (p)$ ; 

- 2) pour toute repr\'esentation lin\'eaire $\rho$ de $G$,
$\rho \circ \mu$ est 
compatible aux r\'eseaux fortement divisibles $M$ de $D=D(b,\mu ,\rho)$ :

$$ M = \oplus_{i \in \Z} D^i_{\mu} \cap M\ .$$

Le groupe \`a un param\`etre  $\mu$ d\'efinit alors un $\otimes$-scindage 
de la filtration du foncteur fibre module sous-jacent 
de la sous-$\otimes$-cat\'egorie pleine de  $\mathrm{MF}_{O_F}^{\mathrm{a}}$
dont les objets sont les r\'eseaux fortement divisibles des
modules de Dieudonn\'e filtr\'es $D(b,\mu,\rho)$,
$\rho$ repr\'esentation lin\'eaire de dimension finie de $G$
(pour la terminologie des $\otimes$-cat\'egorie, voir  \cite{[Sa72]}).
 On impose de plus, que ce scindage s'\'etende
en un $\otimes$-scindage de la sous-$\otimes$-cat\'egorie
pleine de $\mathrm{MF}_{O_F}^{\mathrm{a}}$
dont les objets sont les sous-quotients de ces r\'eseaux.
En particulier, si $\rho_1$ et $\rho_2$ sont deux 
repr\'esentations lin\'eaires de dimension finie de $G$,
si $M_1$ et $M_2$ sont deux r\'eseaux fortement divisibles
de $D(b,\mu, \rho_1)$ et  $D(b,\mu, \rho_2)$ respectivement,
et $M'_1$ et $M'_2$ deux quotients de $M_1$ et $M_2$
dans la cat\'egorie $\mathrm{MF}_{O_F}^{\mathrm{a}}$,
les noyaux de $M_1 \rightarrow M'_1$ et $M_2 \rightarrow
M'_2$ sont compatibles aux graduations ; on impose de plus que  si
$f: M'_1 \rightarrow M'_2$ est un morphisme dans 
$\mathrm{MF}_{O_F}^{\mathrm{a}}$, il est
compatible aux graduations induites sur $M'_1$
et $M'_2$ respectivement.

\subsection{}\label{bg}

Soit $B(G)$ l'ensemble $G(L)$ modulo $\sigma$-conjugaison :
$g\sim g'$ s'il existe $h\in G(L )$ tel que 
$g = hg' \sigma (h^{-1} )$ (\cite{[K85]}, \cite{[K97]}).
Soit $\overline{F}$ la cl\^oture alg\'ebrique
de $F$ dans $\overline{L}$ et $\Gamma_F$ le groupe de
Galois de $\overline{F}/F$.  
Le compos\'e de l'homomorphisme $\widetilde{\kappa}_G : G(L)\rightarrow 
\pi _1(G)_{I_L}$ que nous avons d\'ej\`a 
utilis\'e dans la d\'emonstration de la 
proposition \ref{pargalprop} et 
de   $\pi _1(G)_{I_L}    \rightarrow \pi _1(G)_{\Gamma_F}$
est l'invariant  de Kottwitz
$\kappa_G : B(G)\rightarrow \pi _1(G)_{\Gamma_F}$.

Pour  $\mu $ groupes \`a un param\`etre de 
$G_{\overline{L}}$, rappelons la d\'efinition du sous-ensemble
$B(G, \mu )$ de $B(G)$
(\cite{[K97]} et \cite{[RR96]}).
Soit $T$ un tore maximal de $G_{\overline{F}}$.
Soit $\overline{C}\subset X_* (T)$ une chambre de Weyl
ferm\'ee. 
On note $\mu^{*}$ l'unique conjugu\'e de 
$\mu$ qui appartient \`a $\overline{C}$.
Comme le couple $(T,\overline{C})$ est canoniquement 
associ\'e \`a $G$, on dispose d'une action naturelle de 
$\Gamma_F$ sur $X_* (T)$ qui laisse stable $\overline{C}$.
Si 
$\Gamma_{\mu^{*}}$ est le fixateur de $\mu^*$
dans le groupe de Galois $\Gamma_F$, on pose  :

$$  \overline{\mu}^{*}={1\over (\Gamma_F :\Gamma_{\mu^{*}})}
\times \sum_{\tau\in \Gamma_F / \Gamma_{\mu^*}} \tau(\mu^* ).$$

Soit $b\in G(L)$ et $\overline{\nu _b}\in \overline{C}$
son point de Newton. 
La premi\`ere condition d\'efinissant $B(G,\mu )$ est
que $ \overline{\mu}^{*}\geq \overline{\nu _b}$, autrement
dit que  $ \overline{\mu}^* -  \overline{\nu _b}$ est combinaison
lin\'eaire \`a coefficients $\geq 0$ de coracines simples de $G_{\overline{L}}$
relatives \`a $\overline{C}$. Dans le cas du groupe lin\'eaire,
cette condition dit que le polygone de Newton est au dessus
du polygone de Hodge. D'autre part, soit $\mu^{\natural}$
l'image de $\mu^{*}$ dans $\pi_1(G)_{\Gamma_F}$. La seconde
condition d\'efinissant $B(G,\mu )$ est que 
$\kappa _G (b)= \mu^{\natural}$. Si $G=\GL$, elle 
dit que les polygones de Hodge et de Newton 
ont m\^emes extr\'emit\'es. Noter que ces conditions
ne d\'ependent que de la classe de $\sigma$-conjugaison
$[b]$ de $b$ et de la classe de conjugaison $\{ \mu \}$ dans
$G_{\overline{L}}$. On utilisera aussi la notation 
$B(G,\{ \mu \})$.

\begin{prop}\label{superadm} Soit $(b,\mu )$ une paire
admissible avec $\mu$ d\'efini sur $L$. 
On suppose que $[b]\in B(G,\mu)$.
Alors, il existe un groupe \`a un param\`etre $\mu '$,
d\'efinissant la m\^eme filtration que $\mu$, 
et un $g\in G(L)$ tels que, si $(b_1, \mu_1)$
est le conjugu\'e de $(b, \mu ')$ par $g$,
$(b_1, \mu_1)$ soit super-admissible.
\end{prop} 

\subsection{Remarque.}\label{remsuperadm} Supposons
$G$ quasi-d\'eploy\'e sur $F$ et soit $L'$ une extension
finie de $L$ contenue dans $\overline{L}$.
Fontaine et Rapoport ont prouv\'e dans \cite{[FR02]}
que, si $b\in G(L )$ et $\{ \mu \}$ est une classe
de conjugaison de groupes \`a un param\`etre de 
$G$ d\'efinie sur $L'$, il existe $\mu\in \{ \mu \}$ d\'efinie
sur $L'$ tel que 
$(b,\mu )$ soit admissible si et seulement si 
l'on a la premi\`ere condition d\'efinissant 
$B(G,\mu )$ : $ \overline{\mu}\geq \overline{\nu _b}$.
On voit donc que dans la proposition, on pourrait 
remplacer la condition  $[b]\in B(G,\mu)$ par  
$\kappa _G (b)= \mu^{\natural}$ (\it cf \ref{bg}).\rm

\subsubsection{} L'\'enonc\'e suivant r\'esulte imm\'ediatement 
de la proposition et du th\'eor\`eme de Fontaine et Rapoport
rappel\'e dans la remarque ci-dessus. 

\begin{cor}\label{existesuperadm} 
On suppose $G$ quasi-d\'eploy\'e sur $F$.
Soient $\{\mu \}$ une classe
de conjugaison de groupes \`a un param\`etre de $G$ d\'efinie sur $L$ et 
$[b]$ une classe de $\sigma$-conjugaison 
de $G(L )$. On suppose que $[b]\in B(G,\{ \mu\})$. Alors,
il existe $b\in [b]$ et $\mu \in \{\mu\}$ tels 
que $(b,\mu)$ soit super-admissible.
\end{cor}

\subsection{Preuve de la proposition.}
Nous avons prouv\'e dans \cite{[W84]} qu'il existe un 
$\otimes$-scindage de la filtration du foncteur
fibre sous-jacent de $\mathrm{MF}_{O_F}^{\mathrm{a}}$.
Les graduations que ce scindage 
induit sur les $D(b,\mu',\rho)$ d\'efinissent un 
groupe \`a un param\`etre $\mu'$ de $G$, d\'efini sur 
$L$, qui induit la m\^eme filtration que $\mu$ 
et est tel que $(b,\mu ')$ v\'erifie la condition
2) de la d\'efinition de la super-admissibilit\'e.

\begin{lemme}\label{c}  Posons $c= b \sigma (\mu ' (\pi^{-1} ))\in G(L)$.
Alors, $c$ est $\sigma$-conjugu\'e \`a l'identit\'e :
il existe $g\in G(L)$ tel que $c= g^{-1}\sigma (g)$.
\end{lemme}

\emph{D\'emonstration.} Soit  $[c]$ la classe de $c$ dans 
$B(G)$.
Par Kottwitz (4.13 de \cite{[K97]}),
il s'agit de prouver que $\kappa_G ([c])=0$ et que 
le morphisme de Newton $\nu_c$ est trivial.

Prouvons que $\kappa_G ([c])=0$. 
Comme $\widetilde{\kappa}_G$ est un fonctoriel pour les 
morphismes de groupes $G$, et que, pour $\lambda\in L^*$,  
$\widetilde{\kappa}_{\Gm}(\lambda )=v(\lambda )$, 
o\`u $v$ est la valuation de $L$ normalis\'ee par $v(\pi )=1$,
 $\widetilde{\kappa}_G(\mu '(\pi^{-1}))$ est l'image de 
$-\mu '$ dans $\pi_1 (G)_{I_L}$. Comme $\mu$ et 
$\mu'$ sont conjugu\'es par un automorphisme int\'erieur
de $G$, c'est aussi l'image $-\widetilde{\mu}$ de $-\mu$. 
Comme $\widetilde{\kappa}_G$
est un homomorphisme de groupes, on a donc :

$$\widetilde{\kappa}_G (c )=\widetilde{\kappa}_G (b)-\sigma 
(\widetilde{\mu} ).$$

Comme $\sigma$ agit trivialement sur $\pi_1 (G)_{\Gamma_F}$, on 
voit donc que l'image $\kappa_G ([c])$ de $\widetilde{\kappa}_G (c )$
dans $\pi_1 (G)_{\Gamma_F}$ est $\kappa (b)-\mu^{\natural}$, soit $0$
puisque $[b]\in B(G,\{\mu \})$.            

Prouvons que le morphisme de Newton $\nu_c$ est trivial.
Soient $\rho_U$ une repr\'esentation
lin\'eaire de $G$ dans un $F$-espace
vectoriel $U$ de dimension finie, et posons $D=D(b,\mu ,\rho)$. 
Comme $D$ est admissible, un th\'eor\`eme de Laffaille 
dit que $D$ poss\`ede un r\'eseau fortement divisible $M$
(\cite{[L80]}). 

\begin{lemme}\label{fd} Pour tout r\'eseau fortement divisible
$M$ de $D$, on a : $M=\rho (c)\sigma (M)$.
\end{lemme}

\emph{D\'emonstration du lemme \ref{fd}.}
On a, si $M^i= D_{\mu '}^i \cap M$,
$M=\oplus_{i\in Z} M^i$ puisque $M$ est compatible \`a la graduation
de $D$ d\'efinie par $\mu '$. Puisque $M$ est fortement divisible, on a,
$\mathcal{F}_M ^i = \sum_{i'\geq i} M^{i'}$  :

$$M=\phi (\sum_{i\in\Z}\pi^{-i}\mathcal{F}_M ^i )=\phi 
(\sum_{i\in\Z}\pi^{-i} M ^i )=\phi \mu '(\pi^{-1})(M).$$

Puisque  $\phi=\rho (b)(\sigma\otimes \mathrm{id})$,
il en r\'esulte que  :

$$ M=
\rho(b)(\sigma\otimes \mathrm{id})\mu '(\pi^{-1})(M).$$  

Comme $c= b \sigma (\mu ' (\pi^{-1} ))$, on voit que l'on a bien :
$M=\rho (c)\sigma (M)$.

Il en r\'esulte 
que l'isocristal $(L\otimes_F U , \rho (c)(\sigma\otimes \mathrm{id}))$
est de pente $0$. Comme c'est vrai pour tout $\rho$, 
$\nu _c$ est trivial (3 et 4 de \cite{[K85]})
, ce qui ach\`eve la d\'emonstration du  lemme \ref{c}.

Prouvons la proposition. Soient $b_1=gb\sigma(g^{-1})$
et $\mu_1=\mathrm{int}(g)(\mu ' )$. Posons $c_1 =b_1 \sigma (\mu_1 
(\pi^{-1}))$.
On a :

$$c_1= gb\sigma (g^{-1})\sigma (g)\sigma (\mu' (p^{-1}))\sigma (g^{-1})
=gc\sigma (g^{-1}) =gg^{-1}\sigma (g)\sigma (g^{-1})=1.$$

Le couple $(b_1,\mu_1)$ v\'erifie donc la condition 1) de la 
d\'efinition de la super-admissibilit\'e.
Comme $(b,\mu ')$ v\'erifie la condition 1), la proposition 
r\'esulte alors du fait que $\rho(g)$ induit
un isomorphisme de $D(b,\mu ',\rho)$ sur  $D(b_1,\mu _1,\rho)$
(proposition \ref{conj}).

\begin{cor}\label{mt} Soit $(b,\mu)$ super-admissible.
Pour toute repr\'esentation lin\'eaire $\rho_U$ de $G$
dans un $F$-espace vectoriel de dimension finie $U$,
les r\'eseaux fortement divisibles $M$  de $D(b,\mu,\rho )$ sont les 
r\'eseaux  qui v\'erifient les deux propri\'et\'es suivantes :

a) $M$ provient par extension des scalaires de 
$O_F$ \`a $O_L$ d'un r\'eseau de $U$ ;

b) $M$ est compatible \`a la graduation d\'efinie par 
$\mu$ : $M=\sum_{i\in \Z} M\cap D_{\mu}^i$.
\end{cor}

\emph{Preuve.} Supposons $(b,\mu)$ super-admissible.
La condition a)  r\'esulte imm\'ediatement 
du lemme \ref{fd}, puisque, comme $(b,\mu)$ est super-admissible,
$c=1$. La condition b) fait partie de la d\'efinition de 
la super-admissibilit\'e.

R\'eciproquement, si $M$ est un r\'eseau v\'erifiant a), on
$\sigma (M)=M$. Si  $M$ v\'erifie b), on a :

$$\sum_{i\in \Z} \pi^{-i} (M\cap \mathcal{F}^i _D) = 
 \mu(\pi^{-1})(M).$$

Si $M$ v\'erifie a) et b), on a donc :

$$\phi (\sum_{i\in \Z} \pi^{-i} (M\cap \mathcal{F}^i _D))=
(\sigma\otimes\mathrm{id})\mu (\pi) \mu(\pi^{-1})(M)=\sigma (M)=M,$$

ce qui prouve le corollaire.

\subsection{} Dans la proposition qui suit, on d\'esigne de la m\^eme 
fa\c{c}on un parahorique $K$ de $G$ (resp. $G_L$)
et le sch\'ema en groupes (connexe) sur $O_F$
(resp. $O_L$) dont il est le groupe des points.
Si $K$ est un parahorique de $G$, on d\'esigne par 
$\overline{K}$ la r\'eduction modulo $\pi$ du $O_F$-sch\'ema 
en groupe $K$, et $\overline{K}_{\mathrm{red}}$ le quotient de 
$\overline{K}$ par son radical unipotent. 

\begin{prop}\label{K} Soit $(b,\mu)$ super-admissible.
Alors, il existe un parahorique $K$ de $G$
tel que :

- le groupe \`a un param\`etre $\mu$ se prolonge 
en un groupe \`a un param\`etre de $K_{O_L}$ que l'on note
encore $\mu$  ; on d\'esigne par $\overline{\mu}$ sa r\'eduction modulo $\pi$ ;

- soit  $\overline{\mu}_{\mathrm{red}}$ l'image de $\overline{\mu}$
dans $\overline{K}_{\mathrm{red}}$ ; alors,
 les groupes \`a un param\`etre 
$\sigma^i (\overline{\mu}_{\mathrm{red}})$, pour 
$i\in\Z$  commutent 
entre eux,  autrement dit l'image de  $\overline{\mu}_{\mathrm{red}}$
est contenue dans un tore de $\overline{K}_{\mathrm{red}}$ d\'efini sur $k_F$.

\end{prop}

\subsection{Remarque.} Soient $K$ un parahorique de $G$ et 
$\mu$ un groupe \`a un param\`etre de $G_L$ qui se prolonge
en un groupe \`a un param\`etre de $K_{O_L}$. 
Soit  $\rho_N: K\rightarrow \GL_N$ une repr\'esentation
lin\'eaire de de $K$ dans un $O_F$-module
libre $N$ de rang fini. Posons $M=O_L\otimes_{O_F} N$,
$b=\sigma (\mu) (\pi)$, et soit $\rho$ la repr\'esentation
lin\'eaire de $G$ obtenue \`a partir de $\rho_N$ par extension des 
scalaires \`a $F$. Alors, il n'est pas difficile de 
prouver que $M$ est un r\'eseau fortement divisible
de $D(b,\mu , \rho )$ et que $(b,\mu )$ est admissible.

\emph{D\'emonstration de la proposition.} 
On note $\mathcal{U}_{\mu }$ le sous-groupe de 
$G(L )$ engendr\'e par les $\sigma ^i (\mu (U_L ))$, $i\in \Z$,
$U_L$ \'etant les unit\'es de $L$. 

\begin{lemme}\label{umu} Le groupe  $\mathcal{U}_{\mu }$  est stable par 
$\sigma$,
born\'e pour la bornologie de $G( L )$ 
(2.2.1 de \cite{[T79]}, 4.2.19 de \cite{[BT72]})
et contenu dans le noyau
de $\widetilde{\kappa} _G$. 
\end{lemme}

\it D\'emonstration.\rm La premi\`ere assertion est claire.
Soit $\rho_U$ une repr\'esentation lin\'eaire fid\`ele
de $G$. D'apr\`es Laffaille (\cite{[L80]}),
le module de Dieudonn\'e filtr\'e $D(b,\mu, \rho)$
poss\`ede un r\'eseau fortement divisible $M$.
Comme $M=\sigma(M)$ et que $M$ est compatible avec la graduation
d\'efinie par $\mu$, $M$ est aussi compatible aux graduations
d\'efinies par les $\sigma^i (\mu )$, $i\in \Z$.  Il en r\'esulte
que $M$ est stable par $\mathcal{U}_{\mu }$. Donc, 
$\mathcal{U}_{\mu }$ est born\'e. Comme $U_L$ est le noyau
de $\widetilde{\kappa}_{\Gm}$, il r\'esulte de la fonctorialit\'e 
de $\widetilde{\kappa}$ (7.4. de \cite{[K97]})  que  
$\mathcal{U}_{\mu }$ est contenu dans le noyau de $\widetilde{\kappa} _G$,
ce qui ach\`eve de prouver le lemme.

\begin{lemme}\label{parahoric}
Le groupe $\mathcal{U}_{\mu }$ stabilise une facette  
de l'immeuble $\mathcal{B}(G_{\mathrm{ad}},L)$
qui est stable par $\sigma$.
\end{lemme}

\it D\'emonstration. \rm Le sous-groupe d'automorphismes
de l'immeuble engendr\'e par $\sigma$ est born\'e
car $\sigma$ a un point fixe dans $\mathcal{B}(G_{\mathrm{ad}},L)$
($\mathcal{B}(G_{\mathrm{ad}},F)$ est non vide,
1.10 de \cite{[T79]}). Comme $\sigma$ normalise  
$\mathcal{U}_{\mu }$ et que  $\mathcal{U}_{\mu }$ est born\'e,
le sous-groupe d'automorphismes de l'immeuble engendr\'e par
$\mathcal{U}_{\mu }$ et $\sigma$ est born\'e. Il a donc
un point fixe $x$. Soit $\mathcal{F}_{x}$ la plus petite facette
contenant $x$. Elle est stable par $\sigma$ et $\mathcal{U}_{\mu }$, ce qui 
prouve le lemme.    

Soit $\mathcal{F} _1$ une facette de $\mathcal{B}(G_{\mathrm{ad}},L)$
comme dans le lemme pr\'ec\'edent. Comme $\mathcal{F} _1$
est stable par $\sigma$, l'intersection $\mathcal{F} _{F,1}$
de $\mathcal{F}_1$ 
avec $\mathcal{B}(G_{\mathrm{ad}},F)$ est une facette de 
$\mathcal{B}(G_{\mathrm{ad}},F)$ (5.1.28 de \cite{[BT84]}).
Soit $K_1$ le parahorique de $G$ d\'efini par cette facette
(5.2.6 de \cite{[BT84]}) : le groupe $K_1 (O_L )$ de ses points
\`a valeurs dans $O_L$ est l'intersection du stabilisateur dans $G(L)$
de $\mathcal{F}_{F,1}$ avec le noyau de $\widetilde{\kappa}_G$
(\cite{[HR?]}, et  voir aussi \ref{G'}).
On voit donc avec les deux lemmes pr\'ec\'edents que 
$\mathcal{U}_{\mu }$ est contenu dans $K_1 (O_L )$. 
Il en r\'esulte que $\mu (U_L )$ est contenu dans 
$K_1 (O_L )$. Comme $\Gm$ est \'etoff\'e sur $O_L$
(1.7 de \cite{[BT84]}), le groupe \`a un param\`etre
$\mu : \Gm\rightarrow G_L$ se prolonge en un groupe
\`a un param\`etre de $(K_{1})_{O_L }$. On note encore
$\mu$ ce prolongement. On a r\'ealis\'e avec $K_1$ la premi\`ere
condition de la proposition.

Notons $\overline{\mu}_1 :\Gm\rightarrow \overline{K_1} $
la r\'eduction modulo $\pi$ de $\mu$ et 
$\overline{\mathcal{U}}_1$ le sous-groupe alg\'ebrique de 
$\overline{K_1}$ engendr\'e par les images des 
$\sigma ^i (\overline{\mu}_1 )$, $i\in \Z$.
C'est un sous-groupe r\'eduit, d\'efini sur
$k_F$ et connexe de $\overline{K_1}$. 
Soit $\rho_1$ une repr\'esentation lin\'eaire fid\`ele
de $K_1$  dans un $O_F$-module
libre de type fini $M_F$ (1.4.5 de \cite{[BT84]}).
Posons $M=O_L \otimes_{O_F} M_F$.
D'apr\`es le corollaire \ref{mt}, 
$M$ est un r\'eseau fortement divisible de 
$D=D(\sigma (\mu (p)),\mu,\rho_1)$. Il d\'efinit
donc un objet de la cat\'egorie $\mathrm{MF}_{O_F}^{\mathrm{a}}$.
Soit $\overline{M}$ l'objet $M/\pi M$. 
Soit $(0)=\overline{M}_0\subset \overline{M}_1\subset
\ldots \overline{M}_r = \overline{M}$ une suite de 
Jordan-H\"{o}lder du module de Dieudonn\'e filtr\'e
$\overline{M}$. Comme les r\'eseaux images r\'eciproques
dans $M$ des $\overline{M}_{r'}$ sont des r\'eseaux fortement
divisibles, il r\'esulte du corollaire \ref{mt} que 
les sous-espaces vectoriels  $\overline{M}_{r'}$  du
$k$-espace vectoriel $\overline{M}$ proviennent par 
extension des scalaires de $k_F$ \`a $k$ de sous-espaces
$\overline{M}_{r',F}$
de $M_F /\pi M_F$. La condition 2) de la d\'efinition
de la super-admissibilit\'e entra\^{\i}ne que $\overline{\mu}_1$
laisse stable les  $\overline{M}_{r'}$. On voit donc 
que $\overline{\mathcal{U}}_1$ laisse stable les $\overline{M}_{r',F}$.
D'apr\`es un th\'eor\`eme de Fontaine et Laffaille
(prop. 4.4 de \cite{[FL83]} ; 2.2 de \cite{[W84]}), un  
objet simple $N$ de  $\mathrm{MF}_{O_F}^{\mathrm{a}}$ est 
\'el\'ementaire, en particulier son alg\`ebre d'endomorphismes
$\mathrm{End}(N)$ est une  extension finie  de $k_F$ et $N$ 
est un $\mathrm{End}(N)$-espace vectoriel de dimension $1$.
Soit $N$ l'un 
des quotients $\overline{M}_{r',F}/\overline{M}_{r'-1,F}$.
La condition 2) de la super-admissibilit\'e
entra\^{\i}ne que  l'action de $\overline{\mu}_1$ sur $N$ commute
\`a celle de $\mathrm{End}(N)$. 
Comme $N$ est un $\mathrm{End}(N)$-espace vectoriel de dimension $1$,
il en r\'esulte que 
l'image de $\overline{\mathcal{U}}_1$ est un sous-tore du groupe
multiplicatif de $\mathrm{End}(N)$. On voit donc que l'image 
de $\overline{\mathcal{U}}_1$ dans le groupe lin\'eaire de
$\oplus_{r'} \overline{M}_{r',F}/\overline{M}_{r'-1,F}$ est 
un tore. Il en r\'esulte que $\overline{\mathcal{U}}_1$
est un groupe r\'esoluble. Notons $\overline{U}$ son radical unipotent.
Il est connexe car $\overline{\mathcal{U}}_1$ l'est. 
D'apr\`es un th\'eor\`eme de Borel et Tits (\cite{[BT71]}),
il existe un sous-groupe parabolique $\overline{P}$ de 
$\overline{K_1}$
qui contient le normalisateur de $\overline{U}$ et dont le radical
unipotent contient $\overline{U}$. L'image r\'eciproque de $\overline{P}(k)$
dans $K_1 (O_L)$ est le groupe des points \`a valeurs
dans $O_L$ d'un parahorique que l'on note $K$ (4.6.33 de 
\cite{[BT84]}). 

Prouvons que $K$ satisfait aux conditions de la proposition.
Comme $\overline{\mathcal{U}}_1$ normalise son radical
unipotent $\overline{U}$, le parabolique $\overline{P}$ contient 
$\overline{\mathcal{U}}_1$
et donc l'image de $\overline{\mu}_1$. Il en r\'esulte
que $K(O_L )$ contient bien l'image de $\mu$. 
Comme l'image de $\overline{U}$ dans 
$\overline{K}_{\mathrm{red}}$ est triviale,
l'image de  $\overline{\mathcal{U}}_1$
dans $\overline{K}_{\mathrm{red}}$ est un tore qui contient
l'image de $\overline{\mu}_1$. Ceci 
ach\`eve de prouver la proposition.
  
\section{Positions relatives de sous-groupes d'Iwahori et de 
leurs images par le Frobenius.}
\label{posiwahori}

\subsection{Rappels (4 de \cite{[R02]}).}\label{rappels}

Dans ces rappels, $G$ est un groupe r\'eductif quelconque
sur $L$.

\subsubsection{La ``composante neutre'' $G(L)'$ de $G(L)$ 
(5.2.11. de \cite{[BT84]}
et \cite{[R02]}).}\label{G'}

Dans ce num\'ero, nous esquissons les preuves de faits
connus (\cite{[HR?]}).

Soit $G(L)'$ le sous-groupe de $G(L)$ engendr\'e par les 
sous-groupes parahoriques de $G(L)$.
Le sous-groupe $G(L)'$ est le noyau de $\widetilde{\kappa}_G
: G(L)\rightarrow \pi_{1}(G)_{I_{L}}$.
Soit, en effet, $\mathcal{A}$ un appartement de l'immeuble
$\mathcal{B}(G_{\mathrm{ad}},L)$. 
Soit $S$ le tore maximal d\'eploy\'e de $G$ 
d\'efini par $\mathcal{A}$ et 
$T$ son centralisateur. Comme le corps r\'esiduel
de $L$ est alg\'ebriquement clos, $G$ est quasi-d\'eploy\'e et
$T$ est un tore maximal de $G$. Le groupe $\pi_{1}(G)$
s'identifie donc au quotient de $X_{*} (T)$ par le 
r\'eseau engendr\'e par les coracines.
Si $I_L$ est le groupe de Galois
de $\overline{L}/L$, soit $T(L)_1$ le noyau de 
$\widetilde{\kappa}_T : T(L)\rightarrow X_{*} (T)_{I_L}$.
Alors, $G(L)'$ est engendr\'e par $T(L)_1$ et des sous-groupes
unipotents (\cite{[BT84]} \it loc. cit.\rm). 
Il en r\'esulte que $G(L)'$ est contenu dans le 
noyau de $\widetilde{\kappa}_G$. Soit $g\in G(L)$. Prouvons que si  
$\widetilde{\kappa}_G (g)=1$, $g\in G(L)'$. Soit $N(L)$ le 
normalisateur de $T$ et $\widetilde{W}= N(L)/T(L)_1$ le groupe
d'Iwahori-Weyl (relatif \`a $S$). Si $\mathrm{Iw}$ est 
un sous-groupe d'Iwahori associ\'e \`a une alc\^ove
contenue dans $\mathcal{A}$,
il existe $n\in N(L)$, $k_1 ,k_2 \in \mathrm{Iw}$ tels 
que $g=k_1 n k_2$. On a 
$\widetilde{\kappa}_G (k_1 )=\widetilde{\kappa}_G (k_2 )=1$
comme on vient de le voir. On a donc 
$\widetilde{\kappa}_G (n )=\widetilde{\kappa}_G (g )=1$.
Soit $W_{\mathrm{a}}\subset \widetilde{W}$ le groupe de Weyl affine
(associ\'e \`a $S$).
Si $G_{\mathrm{sc}}$ est le rev\^etement universel
du sous-groupe d\'eriv\'e de $G$, et $T _{\mathrm{sc}}$ l'image
inverse de $T$ dans $G_{\mathrm{sc}}$, on a avec des notations
\'evidentes : $W_{\mathrm{a}}= N_{\mathrm{sc}}(L)/T_{\mathrm{sc}}(L)_1$.
L'homomorphisme  $\widetilde{\kappa}_{G_{\mathrm{sc}}}$ est trivial
puisque $\pi_1 (G_{\mathrm{sc}})$ l'est. Il en r\'esulte
que la restriction de $\widetilde{\kappa}_{G}$ \`a
$N(L)$ se factorise \`a travers 
$N(L)\rightarrow \widetilde{W}\rightarrow 
\widetilde{W}/W_{\mathrm{a}}$. Ce dernier groupe
s'identifie \`a $\pi_{1}(G)_{I_{L}}$. De plus, l'homomorphisme
$T(L)\rightarrow \widetilde{W}\rightarrow \pi_{1}(G)_{I_{L}}$
est induit par $\widetilde{\kappa}_{T}$ et donc par 
$\widetilde{\kappa}_{G}$.
On voit donc que 
$\widetilde{\kappa}_G (n )=1$ entra\^{\i}ne qu'il existe
$n'$ dans l'image de $N_{\mathrm{sc}}(L)$
et $t\in T(L)_1$ tels que $n=n't$. Puisque $G(L)'$ 
contient l'image de $G_{\mathrm{sc}}(L)$ et $T(L)_1$ 
(5.2.11. de \cite{[BT84]}),
on a bien $n\in G(L)'$.

On voit que la restriction de $\widetilde{\kappa}_G$ \`a $N(L)$
induit par passage au quotient 
un morphisme $
\widetilde{\kappa}_{\widetilde{W}}: \widetilde{W}\rightarrow \pi_1 (G)_{I_L}$.
La preuve esquiss\'ee ci-dessus
entra\^{\i}ne que 
si $g\in G(L)$, $g=k_1 n k_2$ avec $k_1 ,k_2 \in \mathrm{Iw}$,
$n\in N(L)$, alors $\widetilde{\kappa}_G (g)$ co\"{\i}ncide avec
l'image de $n$ par l'homomorphisme 
$N(L)\rightarrow \widetilde{W} \rightarrow \pi_1 (G)_{I_L}$.

Pour $\mathcal{F}$ facette
de $\mathcal{B}(G_{\mathrm{ad}},L)$, le 
sous-groupe parahorique de $G(L)$ d\'efini par 
$\mathcal{F}$ est l'intersection du fixateur de 
$\mathcal{F}$ avec $G(L)'$ (1.2 de \cite{[R02]}) ; c'est aussi l'intersection 
du stabilisateur de $\mathcal{F}$ avec $G(L)'$, puisque $G(L)'$
est le groupe sous-jacent \`a un syst\`eme de Tits tel
que les parahoriques soient des paraboliques 
(prop. 5.2.12 de \cite{[BT84]})
et qu'un groupe parabolique d' un syst\`eme de Tits est 
son propre normalisateur.

Soit $\mathcal{C}$ l'alc\^ove qui d\'efinit $\mathrm{Iw}$.
Notons  $\Omega_{\mathcal{C}}$ le stabilisateur de $\mathcal{C}$
dans $\widetilde{W}$.  Le morphisme $\widetilde{\kappa}_{\widetilde{W}}$
identifie le groupe
$\Omega_{\mathcal{C}}$  \`a $\pi_1 (G)_{I_L}$ et
$\widetilde{W}$  au produit semi-direct 
de $\Omega_{\mathcal{C}}$ avec le sous-groupe distingu\'e $W_{\mathrm{a}}$ 
(\cite{[R02]}).

\subsubsection{L'invariant $\mathrm{inv}(g_1 , g_2)_{\mathcal{C}}$.}
\label{definv}
Soit $( \mathcal{A} ,\mathcal{C} )$ un couple form\'e 
d'un appartement $\mathcal{A}$ et d'une alc\^ove   
$\mathcal{C}$ de $\mathcal{B}(G_{\mathrm{ad}},L)$
contenue  dans $\mathcal{A}$.
On reprend les notations du num\'ero pr\'ec\'edent,
en pr\'ecisant en indice la d\'ependance de 
$\mathcal{A}$ et  $\mathcal{C}$.

Le groupe d'iwahori-Weyl  $\widetilde{W}_{\mathcal{A}}$ ne d\'epend 
que de $\mathcal{A}$, mais pour $( \mathcal{A}' ,\mathcal{C}')$  comme
$( \mathcal{A} ,\mathcal{C} )$, on dispose 
d'un isomorphisme de $\widetilde{W}_{\mathcal{A}}$ sur
$\widetilde{W}_{\mathcal{A'}}$ :
il est d\'efini par  $\mathrm{int}(h)$ pour  $h\in G(L)'$ tel
que $( \mathcal{A}',\mathcal{C}' )=(h \mathcal{A} , h \mathcal{C})$.
Cet isomorphisme ne d\'epend pas du choix de $h$,
car $h$ est bien d\'efini modulo multiplication \`a droite 
par un \'el\'ement $t$ de l'intersection de $\mathrm{Iw}_{\mathcal{C}}$
avec $N_{\mathcal{A}}(L)$ et un tel $t$ agit trivialement sur $\mathcal{A}$
(prop. 2.3.2. de  \cite{[BT72]}). Ces isomorphismes
sont compatibles en un sens \'evident et permettent 
de d\'efinir un groupe $\widetilde{W}$ associ\'e
au groupe $G$, avec, pour chaque  $( \mathcal{A} ,\mathcal{C} )$
un isomorphisme de $\widetilde{W}$ sur $\widetilde{W}_{\mathcal{A}}$.   
Comme $G(L)'$ est le noyau de $\widetilde{\kappa}_G$,
on voit que 
le morphisme $\widetilde{W}_{\mathcal{A}}\rightarrow 
\pi_1 (G)_{I_L}$ d\'efini au num\'ero pr\'ec\'edent d\'efinit 
un morphisme $\widetilde{W}\rightarrow 
\pi_1 (G)_{I_L}$ que l'on note $\widetilde{\kappa}_{\widetilde{W}}$.
Comme au num\'ero pr\'ec\'edent, on note $W_{\mathrm{a}}$ son noyau et 
$\Omega_{\mathcal{C}}$ le stabilisateur de $\mathcal{C}$ dans 
$\widetilde{W}$.

Soient
$g_1$, $g_2$ deux \'el\'ements de $G(L)$. 
La d\'ecomposition de Bruhat d\'efinit une bijection
naturelle de l'ensemble des doubles classes
$\mathrm{Iw}_{\mathcal{C}} g\mathrm{Iw}_{\mathcal{C}}$ dans 
$\widetilde{W}_{\mathcal{A}}=\widetilde{W}$.
L'invariant
$\mathrm{inv}(g_1 , g_2)_{(\mathcal{A} ,\mathcal{C})}$ 
est d\'efini comme \'etant
l'\'el\'ement $w\in \widetilde{W}$ associ\'e \`a 
la double classe de $g_1 ^{-1}g_2$. On a donc, en notant
$\nu_{(\mathcal{A},\mathcal{C})}$ la projection 
$N_{\mathcal{A}} (L)\rightarrow \widetilde{W}$ :

$$\begin{array}{rrcl}
(1)&  \mathrm{inv}(g g_1 ,gg _2)_{(\mathcal{A} ,\mathcal{C})}
& = & \mathrm{inv}(g_1 ,g_2)_{(\mathcal{A} ,\mathcal{C})},    
\mathrm{pour}\ g\in G(L), \\

(2)&  \mathrm{inv}(g_1 k_1 , g_2 k_2)_{(\mathcal{A} ,\mathcal{C})}
& = & \mathrm{inv}(g_1 ,g_2)_{(\mathcal{A} ,\mathcal{C})}, 
\mathrm{pour}\ k_1,\ k_2\in \mathrm{Iw}_{\mathcal{C}}, \\ 

(3)&  \mathrm{inv}(1 , g)_{(\mathcal{A} ,\mathcal{C})} & = & 
\nu_{(\mathcal{A},\mathcal{C})} (g),
\mathrm{pour}\ g\in N_{\mathcal{A}}(L), \\

(4)& \mathrm{inv}(\mathrm{int}(h)(g_1) ,
\mathrm{int}(h)(g_2))_{h(\mathcal{A} ,\mathcal{C})}&
= & \mathrm{inv}(g_1 ,g_2)_{(\mathcal{A} ,\mathcal{C})},
\mathrm{pour}\ h\in G(L)'.
\end{array}$$

Il en r\'esulte que $\mathrm{inv}(g_1 , g_2)_{(\mathcal{A} ,\mathcal{C})}$
ne d\'epend en fait que de $\mathcal{C}$. En effet, si 
$\mathcal{A}'$ est un appartement contenant $\mathcal{C}$, il existe
$h\in \mathrm{Iw}_{\mathcal{C}}$ tel que 
$\mathcal{A}'=h\mathcal{A}$ et,  on a en utilisant (4) :

$$\mathrm{inv}(g_1 , g_2)_{(\mathcal{A} ,\mathcal{C})}=
\mathrm{inv}(\mathrm{int}(h)(g_1 ) , 
\mathrm{int}(h)(g_2)_{(\mathcal{A}' ,\mathcal{C})},$$

et donc, en utilisant (1) et (2) :  

$$\mathrm{inv}(g_1 , g_2)_{(\mathcal{A} ,\mathcal{C})}=
\mathrm{inv}(g_1 , g_2)_{(\mathcal{A}' ,\mathcal{C})}.$$

On utilisera la notation :  $\mathrm{inv}(g_1 , g_2)_{\mathcal{C}}$.
On voit que :

$$ \widetilde{\kappa}_{\widetilde{W}}(\mathrm{inv}(g_1 , g_2)_{\mathcal{C}})=
\widetilde{\kappa} _{G}(g_2 g_1^{-1}).$$

\subsubsection{ Remarque.}\label{remposiw} 
L'invariant $\mathrm{inv}(g_1 , g_2)_{\mathcal{C}}$
mesure la position relative de sous-groupes d'Iwahori. 

Soient en effet $\mathrm{Iw}_1$ et $\mathrm{Iw}_2$
deux sous-groupes d'Iwahori de $G(L )$. On peut d\'efinir
leur position relative $\mathrm{inv}(\mathrm{Iw}_1,\mathrm{Iw}_2)
\in W_{\mathrm{a}}$ de la fa\c{c}on suivante. On choisit
une alc\^ove $\mathcal{C}$.  Si $g_1$ et $g_2$ sont deux \'el\'ements
de $G(L)'$ tels que $ \mathrm{Iw}_1=g_1\mathrm{Iw}_{\mathcal{C}}$ et  
$ \mathrm{Iw}_2=g_2\mathrm{Iw}_{\mathcal{C}}$, on pose :

$$ \mathrm{inv}( \mathrm{Iw}_1,\mathrm{Iw}_2)=
\mathrm{inv}( g_1,g_2)_{\mathcal{C}}.$$

Cette d\'efinition ne d\'epend pas du choix de $g_1$
et $g_2$ car $\mathrm{Iw}_{\mathcal{C}}$ est son propre normalisateur
dans $G(L)'$. On voit avec les (1) et (4) du num\'ero
pr\'ec\'edent qu'elle
ne d\'epend pas du choix de $\mathcal{C}$.

Pour  $g_1\in G(L)'$ et $g_2\in G(L )$,
la  projection de   $\mathrm{inv}( g_1,g_2)_{\mathcal{C}}$
sur $W_{\mathrm{a}}$ dans la bijection 
$\widetilde{W}=W_{\mathrm{a}}\rtimes \Omega_{\mathcal{C}}$
est $\mathrm{inv}(g_1 \mathrm{Iw},g_2 \mathrm{Iw})$.
En effet, soient $\mathcal{A}$ un appartement contenant
$\mathcal{C}$ et $g'_2 \in G(L)'$, $n_2$ un \'el\'ement
de $N_{\mathcal{A}}(L )$ qui stabilise $\mathcal{C}$
tels que $g_2= g' _2 n_2$ et d'autre part $k_1 ,\ k_2,\ n$ tels que
$g_1 ^{-1}g' _2 = k_1 n k_2$, $n\in N_{\mathcal{A}}(L)$, 
$k_1 ,\ k_2\in \mathrm{Iw}_{\mathcal{C}}$. On a : 

$$ g_1 ^{-1} g_2 = k_1 nn_2 (\mathrm{int}(n_2 ^{-1})(k_2 )).$$

Comme $\mathrm{int}(n_2 ^{-1})(k_2 )\in \mathrm{Iw}_{\mathcal{C}}$,
il en r\'esulte que, en notant $\omega\in \Omega_{\mathcal{C}}$ l'image
de $n_2$ dans $\widetilde{W}$ :

$$ \mathrm{inv}( g_1,g_2)_{\mathcal{C}}=\mathrm{inv}( g_1,g'_2)_{\mathcal{C}}
\omega .$$

Comme on a suppos\'e que $g_1 \in G(L)'$, on voit que la projection 
de $\mathrm{inv}( g_1,g_2)_{\mathcal{C}}$ sur $W_{\mathrm{a}}$
est bien :

$$\mathrm{inv}( g_1 \mathrm{Iw}_{\mathcal{C}},g'_2\mathrm{Iw}_{\mathcal{C}})=
\mathrm{inv}( g_1 \mathrm{Iw}_{\mathcal{C}},g_2\mathrm{Iw}_{\mathcal{C}}).$$

\subsubsection{Invariants par rapport \`a un parahorique
quelconque (\cite{[R02]}). }\label{paraquelconque}

Soient $\mathcal{F}$ une facette de $\mathcal{B}(G_{\mathrm{ad}},L)$
et $K_{\mathcal{F}}=K$ le parahorique qui lui est associ\'e.
Si  $\mathcal{A}$ est un appartement contenant $\mathcal{F}$,
le groupe 
de Weyl $\widetilde{W}^{K}_{\mathcal{A}}$ est 
le quotient $(N_{\mathcal{A}}(L)\cap K)/T_{\mathcal{A}}(L )_1$.
Soient  $(\mathcal{A},\mathcal{C})$ 
un couple form\'e d'un appartement $\mathcal{A}$ contenant $\mathcal{F}$
et d'une chambre $\mathcal{C}$ contenue dans $\mathcal{A}$
et dont l'adh\'erence contient $\mathcal{F}$
et  $(\mathcal{A}',\mathcal{C}')$ comme $(\mathcal{A},\mathcal{C})$. 
Il existe $h\in K$ tel que 
$(\mathcal{A}',\mathcal{C}')=h(\mathcal{A},\mathcal{C})$
et $\mathrm{int}(h)$ induit un isomorphisme de  
$\widetilde{W}^{K}_{\mathcal{A}}$ sur 
$\widetilde{W}^{K}_{\mathcal{A}'}$. Ces isomorphismes
sont compatibles en le sens \'evident et permettent 
de d\'efinir le sous-groupe $\widetilde{W}^{K}$ 
 de $\widetilde{W}$. On a les d\'ecompositions 
en doubles classes :

$$K\backslash G(L)/K\simeq
\widetilde{W}^{K}_{\mathcal{A}}\backslash 
\widetilde{W}_{\mathcal{A}}/
\widetilde{W}^{K}_{\mathcal{A}}.$$
 
Pour $g_1,\ g_2 \in G(L )$ elles permettent de d\'efinir 
$\mathrm{inv}(g_1 ,g_2 )_{\mathcal{F}}\in 
\widetilde{W}^{K}\backslash \widetilde{W}/
\widetilde{W}^{K}$ (on utilisera aussi la notation
$\mathrm{inv}(g_1 ,g_2 )_{K}$). Pour toute chambre 
$\mathcal{C}$ dont l'adh\'erence contient 
$\mathcal{F}$, $\mathrm{inv}(g_1 ,g_2 )_{\mathcal{F}}$
est l'image de $\mathrm{inv}(g_1 ,g_2 )_{\mathcal{C}}$
dans $\widetilde{W}^{K}\backslash \widetilde{W}/
\widetilde{W}^{K}$.

\subsubsection{Les vari\'et\'es de Deligne-Lusztig affines.}
\label{DeligneLusztigaff} 

On suppose de nouveau que $G$ est d\'efini sur $F$.
Soit $K$ un parahorique de $G_F$.
Soit $b\in G(L)$. 
Pour $g\in G(L)$,   
$\mathrm{inv}(g ,b\sigma (g))_{K(O_L)}$
ne d\'epend que de la classe de $g$ dans 
$G(L)/ K(O_L )$.
Si $K$ est un Iwahori, cela
r\'esulte du (2) de \ref{definv} 
et pour $K$ quelconque, on le voit 
de m\^eme. 
Pour $\overline{w}\in \widetilde{W}^{K}\backslash \widetilde{W}/
\widetilde{W}^{K}$,
la vari\'et\'e de Deligne-Lusztig affine g\'en\'eralis\'ee 
$X_w (b)_{K}$ est 
d\'efinie par :

$$ X_{\overline{w}} (b)_{K}= \{ \overline{g}\in G(L)/ 
K(O_L );\ 
\mathrm{inv}(g,b\sigma (g))_{K(O_L)}=\overline{w}\},$$

(il est conjectur\'e que l'on peut munir cet ensemble d'une 
structure de vari\'et\'e alg\'ebrique sur le corps r\'esiduel
$k$ de $L$).

Si $b'$ est $\sigma$ conjugu\'e de $b$ \it i.e. \rm
$b'= h b\sigma (h^{-1} )$, l'application $g \mapsto hg$
d\'efinit une bijection de 
$X_{\overline{w}} (b)_{K}$ sur 
$X_{\overline{w}} (b')_{K}$.
En particulier, si $\mathrm{Iw}$ et $\mathrm{Iw}'$
sont Iwahori de $G_F$, ils sont conjugu\'es par 
$\mathrm{int}(h)$ pour un $h$ dans le sous-groupe $G(F)'$
de $G(F)$ engendr\'e par les sous-groupes parahoriques
de $G(F)$. Comme clairement $G(F)'\subset G(L)'$,
il r\'esulte des (1) et  (4) de \ref{definv}
que  $g\mapsto gh^{-1}$ d\'efinit une bijection de 
$X_{\overline{w}} (b)_{\mathrm{Iw}}$ sur $X_{\overline{w}} (b)_{\mathrm{Iw}'}$.

\subsubsection{Remarque.} 
Comme $G$ est r\'esiduellement quasi-d\'eploy\'e (1.10.2 de
\cite{[T79]}), il existe une alc\^ove de $\mathcal{B}(G_{\mathrm{ad}},L)$ 
qui est stable par $\sigma$. Soit $\mathcal{C}$ une telle
alc\^ove.
Supposons $g\in G(L)'$. Notons $\phi$
l'action de $b\sigma$ sur l'immeuble $\mathcal{B}(G_{\mathrm{ad}},L)$
et en particulier l'ensemble des sous-groupes d'Iwahori de $G_L$.
Si $\mathrm{Iw}$ est le sous-groupe d'Iwahori $g \mathrm{Iw}_{\mathcal{C}}$,
il r\'esulte de la remarque \ref{remposiw}
que la projection de $\mathrm{inv}(g ,b\sigma (g))_{\mathcal{C}}\in 
\widetilde{W}$ sur $W_{\mathrm{a}}$ est 
$\mathrm{inv}(\mathrm{Iw},\phi \mathrm{Iw})$. D'autre part,
l'image de $\mathrm{inv}(g ,b\sigma (g))_{\mathcal{C}}$ 
par  $\widetilde{\kappa}_{\widetilde{W}}: \widetilde{W}
\rightarrow \pi_1 (G)_{I_L}$ est $\widetilde{\kappa}_G (b)$.

\subsection{Calculs d'invariants.}

On suppose toujours que $G$ est un groupe r\'eductif 
sur $F$. On se donne une classe de 
conjugaison de groupes \`a un param\`etre $\{ \mu\}$
de $G_{L}$  et
$[b]\in B(G,\{ \mu\} )$.

\subsubsection{}\label{3.2.1}
Soient $S$ un tore maximal d\'eploy\'e de $G_L$,
$T$ le centralisateur de $S$ et $N$ le normalisateur de $T$. 
Notons $W_0 =N(L)/ T(L )$ le groupe de Weyl relatif
et $\Lambda_T ( \{ \mu\} )$ l'intersection de $X_* (T)$
avec l'ensemble des groupes \`a un param\`etre de
$\{ \mu \}$ qui sont d\'efinis sur $L$.
La proposition suivante 
est sans doute bien connue :

\begin{prop}  $\Lambda_T ( \{ \mu\} )$ est une orbite sous 
$W_0$.\end{prop}

\emph{D\'emonstration.} 
Soient $\mu_1$ et $\mu_2$ deux groupes \`a un param\`etre
de $\Lambda_T ( \{ \mu\} )$.
Comme tout espace homog\`ene 
sous $G(L)$ a un point rationnel sur $L$ (cor. 1 du 3.2.3. de 
\cite{[CG]}), il existe $g\in G(L)$ tel que 
$\mu_2 =\mathrm{int}(g)(\mu_1 )$. La proposition 
r\'esulte alors de ce que $S$ et $\mathrm{int}(g)(S)$ sont 
conjugu\'es par un \'el\'ement d\'efini sur $L$ du  centralisateur de $\mu_2$.

Notons $\Lambda ( \{ \mu\} )$ l'image de 
$\Lambda_T ( \{ \mu\} )$ dans $X_* (T)_{I_L}$. Si on choisit une 
alc\^ove dans l'appartement d\'efinie par $S$, 
$\Lambda ( \{ \mu\} )$ s'identifie \`a un sous-ensemble 
de $\widetilde{W}$ (\ref{definv}). Il r\'esulte
de la proposition ci-dessus que $\Lambda ( \{ \mu\} )$
co\"{\i}ncide avec l'ensemble not\'e de fa\c{c}on identique
dans le 3 de \cite{[R02]}. 
Notons, pour $K$ 
parahorique de $G$, $X_{ \{ \mu\}}(b)_K$ la r\'eunion
des $X_{\overline{w}}(b)_K$, pour $\overline{w}$
d\'ecrivant l'image de $\Lambda ( \{ \mu\} )$
dans $W^K \backslash \widetilde{W} / W^K$.

\begin{theo}\label{bonpara} Soit $b\in [b]$. On suppose que 
$G$ est quasi-d\'eploy\'e sur $F$ (et que $\{ \mu\}$
est d\'efinie sur $L$). Alors,
il existe un parahorique $K$ de 
$G_F$
tel que $X_{ \{ \mu\}}(b)_K$ soit non vide.\end{theo}

\emph{D\'emonstration.} Le corollaire \ref{existesuperadm}
entra\^{\i}ne que l'on peut supposer $(b,\mu )$ super-admissible,
avec $\mu \in \{ \mu\}$.
Soit $K$ comme dans la proposition \ref{K}. Comme
$(b,\mu )$ est super-admissible, on a $b=\sigma (\mu (p))$.
De plus $\mu$ se prolonge en un groupe \`a un 
param\`etre de $K_{O_L}$, donc aussi $\sigma(\mu )$.
Il existe  un tore maximal d\'eploy\'e $S$ de $K$, d\'efini
sur $O_L$, et  qui contient $\sigma (\mu )$ (3.5. de \cite{[T79]}).
Il d\'efinit un appartement $\mathcal{A}$ de 
$\mathcal{B}(G_{\mathrm{ad}},L)$ qui contient 
la facette $\mathcal{F}$ d\'efinie par $K$. Soit 
$\mathcal{C}$ une alc\^ove contenue dans $\mathcal{A}$
dont l'adh\'erence contient $\mathcal{F}$.
Il r\'esulte de (3) de \ref{definv} que 
$\mathrm{inv}(1,b)_{\mathcal{C}}= t_{\sigma (\mu)}$,
$t_{\sigma (\mu)}$ d\'esignant l'image de  $\sigma (\mu )$ dans
$X_* (T)_{I_L}$, $T$ d\'esignant le centralisateur de $S$.
Il en r\'esulte que $\mathrm{inv}(1,b)_K$ est l'image de 
$t_{\sigma (\mu)}$ dans $W^K \backslash \widetilde{W} / W^K$.
On voit donc que $X_{  \sigma (\{\mu\})}(b)_K$ est 
non vide. 

Il  r\'esulte de ce que l'on vient 
de prouver appliqu\'e \`a 
$\sigma^{-1}(\{\mu \})$ et $\sigma^{-1}(b)$
que $X_{  \{\mu\}}(\sigma^{-1}(b))_K$ est non
vide. Comme $b$ et $\sigma^{-1}(b)$ sont $\sigma$-conjugu\'es,
on a	une bijection de  
$X_{  \{\mu\}}(\sigma^{-1}(b))_K$ sur $X_{  \{\mu\}}(b)_K$
(\ref{DeligneLusztigaff}) et  $X_{  \{\mu\}}(b)_K$ est non vide. 
La proposition est prouv\'ee.

\subsubsection{Remarque.} Il n'existe pas toujours
d'Iwahori qui satisfasse la conclusion de la proposition 
comme le montre d\'ej\`a le cas des isocristaux associ\'es 
aux courbes elliptiques supersinguli\`eres (1 de \cite{[R02]}). 
L'exemple de la remarque \ref{12} prouve que, la filtration admissible 
\'etant donn\'ee, on ne peut 
pas toujours trouver un r\'eseau fortement divisible qui 
d\'efinisse un parahorique $K$ qui soit hypersp\'ecial. 

\subsubsection{}\label{g}
Soient $b\in [b]$ et $\mu \in \{ \mu \}$ tels 
que $(b,\mu )$ soit super-admissible. Soit $K$ un parahorique
de $G$ comme dans la proposition \ref{K}. Soit 
$\overline{T}_{\mu}$ le sous-tore de $\overline{K}_{\mathrm{red}}$
engendr\'e par les $\sigma^i (\overline{\mu}_{\mathrm{red}})$
pour $i\in \Z$ (prop.  \ref{K}). 
Soit $C_{\overline{T}_{\mu}}$ le centralisateur de 
$\overline{\mu}_{\mathrm{red}}$ dans $\overline{K}_{\mathrm{red}}$.
C'est un groupe r\'eductif connexe et il est quasi-d\'eploy\'e
car $k_F$ est fini. Soit $\overline{T}_{\mathrm{red}}$
un tore maximal de $C_{\overline{T}_{\mu}}$ qui est contenu
dans un Borel de $C_{\overline{T}_{\mu}}$ ; il contient
$\overline{T}_{\mu}$ puisque ce dernier est contenu dans 
le centre de $C_{\overline{T}_{\mu}}$. Soit $\overline{T}$
un sous-tore de $\overline{K}$ dont l'image dans $\overline{K}_{\mathrm{red}}$
est $\overline{T}_{\mathrm{red}}$. D'apr\`es le th\'eor\`eme
4.1. de l'expos\'e 11 de \cite{[SGA3]}, on peut relever
$\overline{T}$ en un sous-tore $S$ de $K$. 
Comme $\overline{T}$ est un tore maximal de $K$, $\overline{T}_k$
est un tore maximal d\'eploy\'e de $K_k$ et 
$S_L$ est un tore maximal d\'eploy\'e de $G_L$ (3.5. de \cite{[T79]}).
Il d\'efinit
un appartement $\mathcal{A}_1$ de $\mathcal{B}(G_{\mathrm{ad}},L)$.
De plus, le Frobenius $\sigma$ agit sur $\mathcal{A}_1$,
de fa\c{c}on compatible \`a l'action sur le groupe
des cocaract\`eres 
$X_* (S_L)$ du tore $S$. Il d\'efinit
donc un \'el\'ement $\theta_1$ du normalisateur
du groupe de Weyl affine agissant sur $\mathcal{A}_1$.  

Soit $\mathcal{F}$ la facette d\'efinie par $K$. Elle
est contenue dans $\mathcal{A}_1$ et, comme $K$ est d\'efini
sur $O_F$, elle est stable sous l'action de  $\theta_1$.
Soit $\mathcal{C}_1$ une alc\^ove de $\mathcal{B}(G_{\mathrm{ad}},L)$,
qui est contenue dans $\mathcal{A}_1$ et dans l'\'etoile de
$\mathcal{F}$ \it i.e. \rm
dont l'adh\'erence contient $\mathcal{F}$ (4.6.33 de \cite{[BT84]}).
L'alc\^ove $\sigma(\mathcal{C}_1 )$ est aussi
contenue dans $\mathcal{A}_1$ et dans l'\'etoile de $\mathcal{F}$.
Soit $W_{\mathcal{F}}$ le groupe des transformations 
affines de $\R \otimes X_* (S)$ qui est engendr\'e par les 
r\'eflexions par rapport aux murs de $\mathcal{A}_1$
qui contiennent $\mathcal{F}$.
Il s'identifie au groupe de Weyl absolu 
de $\overline{K}_{\mathrm{red}}$ relativement \`a 
$\overline{T}_{\mathrm{red}}$ (3.5.1 de \cite{[T79]}). 
Soit  $w_{\mathcal{C}_1}$ l'\'el\'ement de $W_{\mathcal{F}}$ qui est tel que 

$$w_{\mathcal{C}_1}(\mathcal{C}_1 )=\theta_1 (\mathcal{C}_1 ).$$

Soit $\overline{M}$
un sous-groupe de Levi de $\overline{K}$ qui contient $\overline{T}$.  
Soit $\overline{n}$ un \'el\'ement de $\overline{M}(k)$, qui normalise
$\overline{T}$ et dont l'image dans $W_{\mathcal{F}}$ est 
$w_{\mathcal{C}_1}$. 
Comme le normalisateur $N$ de $S$ dans $K$ est lisse 
sur $O_{F}$ (cor. 5.3 bis de l'expos\'e 11 de \cite{[SGA3]}),
il existe 
$n\in N(O_{L})$  qui rel\`eve $\overline{n}$. 
Comme $K$ est connexe et lisse sur $O_F$,
il r\'esulte d'un th\'eor\`eme de Lang qu'il existe $g\in K(O_L)$
tel que $n=\sigma (g)g^{-1}$. Soit $\mathcal{C}$ l'alc\^ove
$g^{-1}\mathcal{C}_1$.

\begin{prop}\label{calculinv} L'alc\^ove  $\mathcal{C}$
est stable par $\sigma$. 
Notons $T$ le centralisateur de $S$ et 
$\mu_T$ l'unique
sous-groupe \`a un param\`etre de $S\subset T$ dont l'image 
dans $\overline{K}_{\mathrm{red}}$ est $\mu$. Alors $\mu _T$
et $\mu$ sont conjugu\'es dans $G_L$. On a :

$$\mathrm{inv}(g,b\sigma (g))_{\mathcal{C}}= t_{\sigma (\mu_T )}  
w_{\mathcal{C}_1} ,$$

$t_{\sigma (\mu_T )}$ \'etant l'image de $\sigma (\mu_T )$
dans $X_* (T)_{I_L}$ (on identifie 
$\widetilde{W}$ \`a un quotient du normalisateur 
de $S_{L}$ dans $G_{L}$ gr\^ace \`a l'isomorphisme d\'efini
par $(\mathcal{A}_1,\mathcal{C}_1)$ (\ref{definv})).
\end{prop}

\subsection{Remarque.} Comme $g\mathcal{C}=\mathcal{C}_1$,
on voit que la projection de $\mathrm{inv}(g,b\sigma (g))_{\mathcal{C}}$
sur $W_{\mathrm{a}}$ est 
$\mathrm{inv}(\mathrm{Iw}_{\mathcal{C}_1},\phi \mathrm{Iw}_{\mathcal{C}_1})$.
Comme $g\in G(L)'$, la projection sur $\pi_1 (G)_{I_L}$ est \'egale
\`a l'image de $\sigma (\mu )$.

\emph{D\'emonstration de la proposition.} On a :

$$\sigma (\mathcal{C} )=\sigma (g^{-1})\sigma (\mathcal{C}_1 )
=g^{-1}n^{-1}w_{\mathcal{C}_1} \mathcal{C}_1= 
g^{-1} \mathcal{C}_1 =\mathcal{C},$$

puisque $n$ agit sur $\mathcal{A}_1$ comme $w_{\mathcal{C}_1}$.

Soit $\mathrm{rad_u} (K)$ l'image r\'eciproque 
dans $K(O_L )$ du radical unipotent de $\overline{K}$.
Les groupes \`a un param\`etre $\mu$ et $\mu _T$ ayant 
m\^eme image dans $\overline{K}_{\mathrm{red}}$, il existe
$u\in \mathrm{rad_u} (K)$ tel que $\mu =\mathrm{int}(u)(\mu _T )$
(th. 10.6 de \cite{[LAG]} et 3. de exp. 9 de \cite{[SGA3]}).

D'apr\`es (4) de \ref{definv}, on a, par conjugaison par $g^{-1}$ qui 
appartient \`a $K(O_L)\subset G(L)'$ :

$$\mathrm{inv}(g,b\sigma (g))_{\mathcal{C}}=
\mathrm{inv}(g,gb\sigma (g)g^{-1})_{\mathcal{C}_1}.$$

En multipliant \`a gauche par $g^{-1}$ ((1) de \ref{definv}) et en utilisant 
$\sigma (g)g^{-1}=n$ et $b=\sigma (\mu)(\pi)$ :

$$ =\mathrm{inv}(1,\sigma (\mu)(\pi)n)_{\mathcal{C}_1}.$$

En utilisant $\mu=\mathrm{int}(u)(\mu _T)$ et
$u\in \mathrm{rad_u}(K)(O_L)$, on obtient, en posant 
$u'=\sigma (u )$) :

$$ = \mathrm{inv}(1,u'\sigma(\mu _T)(\pi){u'}^{-1}n)_{\mathcal{C}_1},$$

soit :

$$ = \mathrm{inv}({u'}^{-1}, 
\sigma(\mu _T)(\pi)n(n^{-1}{u'}^{-1}n))_{\mathcal{C}_1},$$

soit, puisque $u'$ et $n^{-1}u'n$ sont des \'el\'ements de
$\mathrm{rad_u}(K)(O_L)\subset \mathrm{Iw}_{\mathcal{C}_1}(O_L)$ :

$$= \mathrm{inv}(1, \sigma(\mu _T)(\pi)n)_{\mathcal{C}_1}.$$

Comme l'image de $n$ dans $W_{\mathcal{F}}$ est $w_{\mathcal{C}_1}$,
ceci prouve la formule \'enonc\'ee.

\section{Triplets $(\mathcal{F},\lambda,\mathcal{C})$ en bonne 
position.}\label{bonneposition}

\subsection{D\'efinition.}\label{defbonneposition}
Dans ce paragraphe, on se donne un syst\`eme de racines r\'eduit
$R$ dans un $\R$-espace vectoriel $V$.  On s'efforce \`a reprendre 
les notations de \cite{[Bbki]}. 
On d\'esigne par $E$ l'espace affine  dont le groupe des  translations est
le dual $V^*$ de $V$ et par $N(W_{\mathrm{a}} )$ le normalisateur 
du groupe de Weyl
affine $W_{\mathrm{a}}$ de $R$ dans le groupe des transformations
affines de $E$.
On se donne  
un \'el\'ement $\theta$ de   $N(W_{\mathrm{a}} )$.
On se 
donne de plus  une facette $\mathcal{F}$ de $E$ (pour 
les murs relativement \`a $W_{\mathrm{a}}$).
On suppose que :

- a) $\theta$ laisse stable $\mathcal{F}$.

On suppose donn\'e $\lambda\in P(R^{\vee})$ un poids 
du syst\`eme dual. 
On note $R_{\mathcal{F}}$ les racines $\vec{\alpha}$
qui sont telles qu'il existe un mur de direction orthogonal 
\`a $\vec{\alpha}$ qui contienne $\mathcal{F}$ ; 
$R_{\mathcal{F}}$ est donc un sous-syst\`eme de racines
de $R$. Soit $R_{\mathcal{F},\lambda,\theta}$ le syst\`eme
de raicines form\'es des racines
de $R_{\mathcal{F}}$ qui sont orthogonales aux $\theta^i (\lambda )$,
pour $i\in \Z$. Il est clair que $\theta$ agit
sur $R_{\mathcal{F}}$ et $R_{\mathcal{F},\lambda,\theta}$. On suppose que :

- b) $\theta$ laisse stable une chambre $C_{\mathcal{F},\lambda,\theta}$
de $R_{\mathcal{F},\lambda,\theta}$.

Se donner une  alc\^ove $\mathcal{C}$ 
dont l'adh\'erence contient $\mathcal{F}$
revient \`a se donner
une chambre $C_{\mathcal{F}}$ de $R_{\mathcal{F}}$, donc un ensemble
de racines positives  de $R_{\mathcal{F}}$.
Si $\mathcal{C}$ est une telle alc\^ove,
on dit que $(\mathcal{F},\lambda,\mathcal{C})$ est \it en 
bonne position \rm si :

-  pour toute racine $\vec{\alpha}>0$ de $R_{\mathcal{F}}$
(relativement \`a $C_{\mathcal{F}}$) telle que 
$\theta^{-1}(\alpha )<0$, on a $\langle \vec{\alpha},\lambda \rangle >0$.

\subsection{Existence d'alc\^oves en bonne position.}
\begin{prop}\label{existebonneposition} 
Soient $R$, $\theta$, $\mathcal{F}$, $\lambda$, 
et $C_{\mathcal{F},\lambda,\theta}$
comme dans le \ref{defbonneposition}.
Alors, il existe une alc\^ove $\mathcal{C}$ 
dont l' adh\'erence
contient $\mathcal{F}$, telle que 
$(\mathcal{F},\lambda,\mathcal{C})$ soit en bonne position, 
et que toute racine positive de $R_{\mathcal{F},\lambda,\theta}$ positive
relativement \`a $C_{\mathcal{F},\lambda,\theta}$ soit positive 
relativement \`a $C_{\mathcal{F}}$.
\end{prop}

\emph{D\'emonstration.} 
Choisir une chambre $C_{\mathcal{F}}$ de $R_{\mathcal{F}}$,
revient \`a choisir un sous-ensemble
$R_{\mathcal{F},+}$ de $R_{\mathcal{F}}$ tel que 
$R_{\mathcal{F},+} \cup - R_{\mathcal{F},+}$ soit une partition de
$R_{\mathcal{F}}$ et 
tel que $\vec{\alpha},\vec{\alpha}'\in R_{\mathcal{F},+}$
et $\vec{\alpha}+ \vec{\alpha}'\in R_{\mathcal{F}}$ entra\^{\i}ne
$\vec{\alpha}+ \vec{\alpha}'\in R_{\mathcal{F},+}$  
(cor. 1 de prop. 20 de 1.7. de chap. 6 \cite{[Bbki]}).
On fait ce choix de la mani\`ere suivante :

- si pour tout $i\in \Z$ $\langle \vec{\alpha},\theta^i(\lambda) \rangle 
=0$, autrement dit si $\vec{\alpha}$ est une racine de 
$R_{\mathcal{F},\lambda,\theta}$,
$\vec{\alpha}\in R_{\mathcal{F},+}$  si elle est positive 
relativement \`a  $C_{\mathcal{F},\lambda,\theta}$ ;

- s'il existe $i$ tel que $\langle \vec{\alpha},\theta^i(\lambda) \rangle 
\not=0$, soit $i_{0}$ le plus petit entier $\geq 0$ pour lequel
ceci se produit,  $\vec{\alpha}\in R_{\mathcal{F},+} $  si 
$\langle \vec{\alpha},\theta^{i_{0}}(\lambda) \rangle 
>0$. 
 
Il est clair que $R_{\mathcal{F},+}$   
d\'efinit bien un ensemble de racines positives
qui contient les racines positives de   $R_{\mathcal{F},\lambda,\theta}$.

Prouvons qu'alors $(\mathcal{F},\lambda,\mathcal{C})$ est en bonne position.
Soit  
$\vec{\alpha}\in R_{\mathcal{F},+}$. On a : 
$\langle \vec{\alpha},\lambda \rangle \geq 0$. 
Il s'agit de prouver que, si $\vec{\alpha}>0$ et 
$\langle \vec{\alpha},\lambda \rangle = 0$, on a 
$\theta^{-1}(\vec{\alpha})>0$. On a  :

- soit $\langle \vec{\alpha},\theta(\lambda) \rangle=0$ ; alors, 

- soit de plus, 
$\vec{\alpha}\in R_{\mathcal{F},\lambda,\theta}$ et 
$\theta^{-1}(\vec{\alpha})>0$ puisque l'ensemble des racines
positives de $R_{\mathcal{F},\lambda,\theta}$ est stable par
$\theta$ ; 

- soit de plus, il existe $i_0\geq 2$ tel que 
$\langle \vec{\alpha},\theta^{i_{0}}(\lambda) \rangle 
>0$ et $\langle \vec{\alpha},\theta^{i}(\lambda) \rangle 
=0$ pour $0\leq i < i_0$ ; alors 
$\langle \theta^{-1}(\vec{\alpha}),\theta^{i_{0}-1}(\lambda) \rangle 
>0$ et  $\langle \theta^{-1}(\vec{\alpha}),\theta^{i}(\lambda) \rangle 
=0 $ pour $0\leq i < i_0-1$ et on a bien $\theta^{-1}(\vec{\alpha})>0$ ;

- soit $\langle \vec{\alpha},\theta(\lambda) \rangle \not= 0$,
donc est $>0$ puisque $\vec{\alpha}$ est $>0$, d' o\`u 
$\langle \theta^{-1}(\vec{\alpha}),\lambda \rangle >0 $ et  
donc $\theta^{-1}(\vec{\alpha})$ est positive.

Dans tous les cas  $\theta^{-1}(\vec{\alpha})$ est positive
et la proposition est prouv\'ee.

\subsection{Bonne position et calcul de longueurs dans $N(W_{\mathrm{a}} )$.}
\label{longueur}

Soient $R$, $\mathcal{F}$, $\lambda$ et $\theta$ comme
dans le \ref{defbonneposition}. 
Soit $\mathcal{C}$ une alc\^ove dont l'adh\'erence 
contient $\mathcal{F}$, de sorte que l'on 
a une partition de $R_{\mathcal{F}}$
en racines positives et n\'egatives . On note $\Omega_{R,\mathcal{C}}$
le stabilisateur dans $N(W_{\mathrm{a}} )$ de 
$\mathcal{C}$, de sorte que $N(W_{\mathrm{a}} )$
s'identifie au produit semi-direct de $\Omega_{R,\mathcal{C}}$
par le sous-groupe distingu\'e  $W_{\mathrm{a}}$.

On note $l_{\mathcal{C}}$ la fonction longueur du groupe de Coxeter 
$W_{\mathrm{a}}$, relativement au  sous-ensemble g\'en\'erateur
form\'e des r\'eflexions par rapport aux murs de $\mathcal{C}$.
On \'etend la fonction $l_{\mathcal{C}}$ au groupe $N(W_{\mathrm{a}} )$
en posant $l_{\mathcal{C}}(\eta )=l_{\mathcal{C}}(w )$
si $\eta =w \gamma$, $w\in W_{\mathrm{a}}$ et 
$\gamma\in \Omega_{R,\mathcal{C}}$. 
L'entier $l_{\mathcal{C}}(\eta )$ est donc le nombre
de murs s\'eparant $\mathcal{C}$ de $\eta (\mathcal {C})$
(lemme 2 du 1.4 du chap. 4 et prop. 17 du 1.6 du chap. 6
de $\cite{[Bbki]}$). 

La proposition suivante est une g\'en\'eralisation de la prop.
1.23 de \cite{[IM65]}, qui est le cas o\`u la facette
$\mathcal{F}$ est r\'eduite \`a un point sp\'ecial.

\begin{prop}\label{iwahorimatsumoto} Notons $t_{\lambda}\in N(W_{\mathrm{a}})$
la translation de vecteur $\lambda$. 
On a :

$$ l_{\mathcal{C}}(t_{\lambda}\theta)=\sum_{\vec{\alpha}}
\mid \langle \vec{\alpha},\lambda \rangle \mid 
+\sum_{\vec{\alpha}}
\mid \langle \vec{\alpha},\lambda \rangle -1 \mid $$

la seconde somme portant
sur l'ensemble $R_{\mathcal{F},\theta}$ des
racines $\vec{\alpha}\in R_{\mathcal{F}}$ 
v\'erifiant $\vec{\alpha}>0$,
$\theta^{-1}(\vec{\alpha} )<0$ ; pour  la premi\`ere, on choisit 
une racine $\vec{\alpha}$ 
dans chaque couple  $\{\vec{\alpha},-\vec{\alpha}\}$
de racines de $R$ tel que ni
$\vec{\alpha}$ ni $-\vec{\alpha}$ n'appartienne \`a $R_{\mathcal{F},\theta}$, 
et on somme sur ces racines $\vec{\alpha}$.  
\end{prop}

\emph{D\'emonstration.} On reprend la d\'emonstration d'Iwahori
et Matsumoto. Soit $x\in \mathcal{C}$ un point proche
de $\mathcal{F}$. La longueur $ l_{\mathcal{C}}(t_{\lambda}w)$
est le nombre de murs rencontrant le segment $[x,\theta (x)+\lambda]$.
Soit $\vec{\alpha}\in R$ et 
$\alpha$ une racine affine de direction
$\vec{\alpha}$ : $\alpha$ est donc une fonction affine
sur $E$ de direction $\vec{\alpha}$ et dont l'hyperplan des
z\'eros est un mur.
Le nombre de murs de direction orthogonale
\`a $\vec{\alpha}$ (ou $-\vec{\alpha}$) est le nombre d'entiers
dans l'intervalle $[\alpha (x), \alpha (\theta(x))+\lambda ]$.
Faisons tendre $x$ vers un point $x_0$ de $\mathcal{F}$.
La longueur du segment tend vers l'entier 
$\mid \langle \vec{\alpha},\lambda \rangle \mid$.
Si $\vec{\alpha}\notin R_{\mathcal{F}}$, $ \alpha (x_0)$
n'est pas un entier, et on voit que la contribution
des murs orhogonaux \`a $\vec{\alpha }$ est bien  
$\mid \langle \vec{\alpha},\lambda \rangle \mid$.
Si $\vec{\alpha}\in R_{\mathcal{F}}$, $ \alpha (x_0)$
est un entier $a$, et l'intervalle 
$[\alpha (x), \alpha (\theta(x))+\lambda ]$ est de la forme
$[a+\epsilon,a+\langle \vec{\alpha},\lambda \rangle +\epsilon ']$
avec $\epsilon$ du signe de $\vec{\alpha}$ et 
$\epsilon '$ de celui de $\theta^{-1}(\vec{\alpha})$.
Pour $\vec{\alpha}\in R_{\mathcal{F}}$ positive,  
dans chacun des quatre cas de signes possibles, le nombre d'entiers
de l'intervalle  $[\alpha (x), \alpha (\theta(x))+\lambda ]$ est :

- si $\langle \vec{\alpha},\lambda \rangle>0, \theta^{-1}(\vec{\alpha})>0$
: $\mid \langle \vec{\alpha},\lambda \rangle \mid$ ;

- si $\langle \vec{\alpha},\lambda \rangle>0, \theta^{-1}(\vec{\alpha})<0$
: $\mid \langle \vec{\alpha},\lambda \rangle \mid -1$ ;

- si $\langle \vec{\alpha},\lambda \rangle<0, \theta^{-1}(\vec{\alpha})>0$
: $\mid \langle \vec{\alpha},\lambda \rangle \mid$ ;

- si $\langle \vec{\alpha},\lambda \rangle<0, \theta^{-1}(\vec{\alpha})<0$
: $\mid \langle \vec{\alpha},\lambda \rangle \mid +1$.

Ceci prouve la proposition.

\subsubsection{Bonne position et ordre de Bruhat.}\label{ordreBruhat} 

On reprend les hypoth\`eses de \ref{defbonneposition}.
Soit $\mathcal{C}$ une alc\^ove dont l'adh\'erence contient $\mathcal{F}$.
On note $\leq_{\mathcal{C}}$ l'ordre partiel sur 
$N(W_{\mathrm{a}} )$ d\'efini par $\eta\leq _{\mathcal{C}} \eta '$
si :

- $\eta$ et $\eta'$ ont m\^eme image dans $\Omega_{R,\mathcal{C}}$ ;

- alors il existe  $w,w' \in W_{\mathrm{a}}$ et
$\gamma\in \Omega_{R,\mathcal{C}}$ tels que 
$\eta=w\gamma$, $\eta'=w'\gamma$ ; la seconde condition 
pour que $\eta\leq _{\mathcal{C}} \eta '$ est que
$w\leq_{\mathcal{C}} w'$ pour l'ordre de Bruhat
sur $W_{\mathrm{a}}$. Pour la d\'efinition 
de l'ordre de Bruhat, voir par exemple \cite{[Bj83]} et sa 
bibliographie :  en particulier, $w\leq_{\mathcal{C}} w'$
si pour une (ou toute) d\'ecomposition r\'eduite de $w'$,
$w' =s_1 s_2 ...s_{l'}$, il existe $1\leq i_1 < i_2 < ...<i_l\leq l'$
tel que $w= s_{i_1}s_{i_2}...s_{i_l}$. 

\begin{prop}\label{bonnepositionbruhat} On suppose       
$(\mathcal{F},\lambda,\mathcal{C})$ en bonne 
position. Alors, on a :

$$l_{\mathcal{C}}(t_{\lambda})=
l_{\mathcal{C}}(t_{\lambda}\theta)+l_{\mathcal{C}}(\theta).$$

Notons $w_{\mathcal{C}}$ l'\'el\'ement de $W_{\mathcal{F}}$
tel que $\theta (\mathcal{C} )=w_{\mathcal{C}}(\mathcal{C})$.
On a : $t_{\lambda}w_{\mathcal{C}}\leq_{\mathcal{C}} t_{\lambda}$.
\end{prop}

\emph{D\'emonstration.} La proposition pr\'ec\'edente donne :

$$l_{\mathcal{C}}(t_{\lambda})=\sum_{\vec{\alpha}}
\mid \langle \vec{\alpha},\lambda \rangle \mid,$$

la somme portant sur un syst\`eme de repr\'esentants 
des classe $\{ \vec{\alpha},- \vec{\alpha}\}$ de $R$, et :

$$l_{\mathcal{C}}(\theta)=\sum_{\vec{\alpha}\in R_{\mathcal{F},\theta}} 1.$$ 

Comme $(\mathcal{F},\lambda,\mathcal{C})$ est en bonne 
position, si $\vec{\alpha}\in R_{\mathcal{F},\theta}$, on a 
$\langle \vec{\alpha},\lambda \rangle >0$, et donc la proposition
pr\'ec\'edente s'\'ecrit :

$$l_{\mathcal{C}}(t_{\lambda}\theta)=\sum_{\vec{\alpha}}
\mid \langle \vec{\alpha},\lambda \rangle \mid 
+\sum_{\vec{\alpha}}
\mid \langle \vec{\alpha},\lambda \rangle  \mid -1.$$

On a bien l'\'egalit\'e de longueurs cherch\'ee. 

Prouvons $t_{\lambda}w_{\mathcal{C}}\leq_{\mathcal{C}} t_{\lambda}$.
Tout d'abord, comme $w_{\mathcal{C}}\in W_{\mathrm{a}}$,
$t_{\lambda}w_{\mathcal{C}}$ et $ t_{\lambda}$ ont m\^eme image
dans $\Omega_{R,\mathcal{C}}$. 

\begin{lemme} Pour $\eta\in N(W_{\mathrm{a}})$ et 
$\gamma\in \Omega_{R,\mathcal{C}}$, on a : 

$$l_{\mathcal{C}}(\eta\gamma)= l_{\mathcal{C}}(\eta)
=l_{\mathcal{C}}(\gamma\eta)$$

et $l_{\mathcal{C}}(\eta)=l_{\mathcal{C}}(\eta^{-1})$.
\end{lemme}

\emph{D\'emonstration.} On a $\eta\gamma(\mathcal{C})=\eta(\mathcal{C})$,
d'o\`u $l_{\mathcal{C}}(\eta\gamma)= l_{\mathcal{C}}(\eta)$.
Comme $\mathrm{int}(\gamma )$ stabilise l'ensemble des 
r\'eflexions par rapport aux murs de $\mathcal{C}$,
on a $l(\mathrm{int}(\gamma )(\eta))=l(\eta )$.
Comme $\gamma\eta =\mathrm{int}(\gamma )(\eta)\gamma$, on a 
bien $l_{\mathcal{C}}(\eta)=l_{\mathcal{C}}(\gamma\eta)$.

Posons $\eta=w\gamma$, $w\in W_{\mathrm{a}}$ et 
$\gamma\in \Omega_{R,\mathcal{C}}$. On a : 
$\eta^{-1}= \mathrm{int}(\gamma^{-1} )(w^{-1})\gamma^{-1}$.
D'o\`u, puisque $l_{\mathcal{C}}(w)=l_{\mathcal{C}}(w^{-1})$
et la premi\`ere partie du lemme,
$l_{\mathcal{C}}(\eta)=l_{\mathcal{C}}(\eta^{-1})$.
Cela prouve le lemme.

Achevons de prouver la proposition. Comme il existe 
$\gamma\in \Omega_{R,\mathcal{C}}$ tel que 
$\theta=w_{\mathcal{C}}\gamma$, on d\'eduit de l'\'egalit\'e
de longueurs d\'ej\`a prouv\'ee :

$$l_{\mathcal{C}}(t_{\lambda})=
l_{\mathcal{C}}(t_{\lambda}w_{\mathcal{C}})+l_{\mathcal{C}}(w_{\mathcal{C}}).$$

Posons $w=w_{\mathcal{C}}^{-1}$ et $\eta=t_{\lambda}w_{\mathcal{C}}$. 
Comme $l_{\mathcal{C}}(w)=l_{\mathcal{C}}(w^{-1})$, l'\'egalit\'e 
ci-dessus devient : 

$$l_{\mathcal{C}}(\eta)+l_{\mathcal{C}}(w)=l_{\mathcal{C}}(\eta w).$$

Posant $\eta=w'\gamma$, $w'\in W_{\mathrm{a}}$ et 
$\gamma\in \Omega_{R,\mathcal{F}}$, on a :
$\eta w= (w'\mathrm{int}( \gamma )(w))\gamma$ et, avec le lemme 
pr\'ec\'edent : 
$l_{\mathcal{C}}(\eta w)=l_{\mathcal{C}}(w'\mathrm{int}( \gamma )(w))$.
L'\'egalit\'e de longueurs ci-dessus s'\'ecrit :

$$l_{\mathcal{C}}(w')+l_{\mathcal{C}}(\mathrm{int}( \gamma )(w))=
l_{\mathcal{C}}(w'\mathrm{int}( \gamma )(w)).$$

Il en r\'esulte que si on multiplie une d\'ecomposition r\'eduite 
de $w'$ par une d\'ecomposition r\'eduite de $\mathrm{int}( \gamma )(w)$,
on obtient une d\'ecomposition r\'eduite de 
$w'\mathrm{int}( \gamma )(w)$. On en d\'eduit que 
$w'\leq_{\mathcal{C}}w'\mathrm{int}( \gamma )(w)$, d'o\`u :
$\eta \leq_{\mathcal{C}} \eta w$ et :

$$t_{\lambda}w_{\mathcal{C}}\leq_{\mathcal{C}}t_{\lambda}.$$

Ceci ach\`eve de prouver la proposition.

\section{Conclusions.}\label{conclusion}

On reprend les notations de \ref{rappels}.
Soit  $\{ \mu \}$ une classe de conjugaison de groupes
\`a un param\`etre de $G_{\overline{L}}$. Rappelons 
(\cite{[KR02]} et \cite{[R02]}) que :

$$ \mathrm{Adm}(\{ \mu\})=\{ w\in\widetilde{W} ; w\leq \lambda
\ \mathrm{pour\ un }\  \lambda\in \Lambda ( \{ \mu\} )\}.$$

Si $\{ \mu \}$ est d\'efinie sur $L$, nous avons donn\'e 
une d\'efinition de 
$\Lambda ( \{ \mu\} )$  au \ref{3.2.1}. L'ordre 
sur $\widetilde{W}$ est d\'efini \`a partir de l'ordre
de Bruhat d'une mani\`ere analogue \`a celle du \ref{ordreBruhat}.
Soit $\mathrm{Iw}$ un sous-groupe
d'Iwahori de $G$.  Soit $b\in G(L)$
tel que $[b]\in B(G,\{ \mu\})$.
Rappelons que
$X(b,\{ \mu\})_{\mathrm{Iw}}$ est la r\'eunion des 
$X_w (b)_{\mathrm{Iw}}$ pour $w\in \mathrm{Adm}(\{ \mu\})$.
Le th\'eor\`eme suivant prouve une conjecture de Kottwitz
et Rapoport sous des hypoth\`eses de non ramification 
(\cite{[KR02]}, \cite{[R02]}).

\begin{theo}\label{conjadm} Supposons $G$ quasi-d\'eploy\'e
et $\{ \mu \}$ d\'efinie sur $L$.  
Alors 
$X(b,\{ \mu\})_{\mathrm{Iw}}$ est non vide. 
Plus pr\'ecis\'ement, il existe $\lambda\in  \Lambda ( \{ \mu\} )$
et $w\in \widetilde{W}$ tels que  $X_w (b)_{\mathrm{Iw}}$
soit non vide et $l(\lambda)=l(w)+l(w^{-1}\lambda)$.
\end{theo}

\subsection{Remarques.}\label{Bruhatfaible}
1) Rappelons que l'ordre de Bruhat faible $\leq '$
sur $W_{\mathrm{a}}$ est d\'efini par $w\leq ' w'$ si
$l(w')=l(w)+l(w^{-1}w')$ (\cite{[Bj83]}).
Si $w\leq ' w'$, on a donc $w\leq  w'$. On voit donc
que l'on prouve en fait la conjecture de Kottwitz et Rapoport
o\`u dans la d\'efinition de  $ \mathrm{Adm}(\{ \mu\})$, 
on remplace l'ordre de Bruhat usuel par l'ordre de Bruhat faible.

2) Si $G$ est non ramifi\'e, il r\'esulte de la 
prop. 4.4 de \cite{[R02]} que, si  $X(b,\{ \mu\})_{\mathrm{Iw}}$
est non vide, $[b]\in B(G,\{ \mu \})$.

\subsection{D\'emonstration du th\'eor\`eme.} 
D'apr\`es le corollaire \ref{existesuperadm},
il existe $b'\in[b]$ et $\mu\in \{ \mu\}$ tels que 
$(b',\mu )$ soit super-admissible. 
On peut donc supposer $(b,\mu)$ super-admissible. 
Soit alors $K$ un parahorique
comme dans la proposition \ref{K}. 
Soit $\mathcal{F}$ la facette 
de $\mathcal{B}(G_{\mathrm{ad}},L)$ qui lui est associ\'ee.
On reprend les d\'efinitions et
notations de \ref{g} et de la proposition \ref{calculinv}.
On a donc  un appartement $\mathcal{A}_1$ de $\mathcal{B}(G_{\mathrm{ad}},L)$
contenant 
$\mathcal{F}$ qui est stable sous l'action 
de $\sigma$ ; l'action de $\sigma$ sur $\mathcal{A}_1$
est not\'ee $\theta_1$. Pour toute alc\^ove 
$\mathcal{C}_1$ contenue dans $\mathcal{A}_1$
dont l'adh\'erence contient $\mathcal{F}$, 
il existe $g\in G(L)$ tel que, si $\mathcal{C}$ est l'alc\^ove
$g^{-1}\mathcal{C}_1$, on ait :   

$$\mathrm{inv}(g,b\sigma (g))_{\mathcal{C}}= t_{\sigma (\mu_T )}  
w_{\mathcal{C}_1},$$

$w_{\mathcal{C}_1}\in W_{\mathcal{F}}$ \'etant d\'efini
par $w_{\mathcal{C}_1}(\mathcal{C}_1 )=\theta_1 (\mathcal{C}_1 ).$
De plus, $\mu$ et $\mu_T$ sont conjugu\'es dans $G_L$, donc 
aussi $\sigma (\mu )$ et $\sigma (\mu_T )$. Il en r\'esulte
que $\sigma (\mu_T )\in \Lambda_{T}(\sigma(\{ \mu \}))$
(\ref{3.2.1}). Pour prouver le th\'eor\`eme
pour l'Iwahori d\'efini par $\mathcal{C}$ et 
la classe de conjugaison de groupes \`a un param\`etre
$\sigma(\{ \mu \})$, il suffit donc de prouver 
que l'on peut choisir $\mathcal{C}_1$ telle que 
$t_{\sigma (\mu_T )}  w_{\mathcal{C}_1}$ soit plus 
petit que $t_{\sigma (\mu_T )}$, pour l'ordre de Bruhat de 
$\widetilde{W}$ identifi\'e au groupe 
d'Iwahori-Weyl relatif \`a $(\mathcal{A}_1,\mathcal{C}_1)$
(\ref{definv}).

Soit $x$ un point sp\'ecial de $\mathcal{A}_1$ de sorte que l'on a 
un syst\`eme de racines r\'eduit $R$ dont 
$W_{\mathrm{a}}$ est le groupe de Weyl affine (1.7. de \cite{[T79]},
prop. 9 du 6.2. de \cite{[Bbki]}).
Comme $\sigma(\mu_T )$ normalise $W_{\mathrm{a}}$,
l'image de $\sigma(\mu_T )$ dans le groupe des transformations
affines de $\mathcal{A}_1$ est la translation par 
un  poids du syst\`eme racines dual de $R$ ; on note ce poids
$\lambda$. On pose : $\theta=\theta_1$ et on reprend
les notations de \ref{defbonneposition}.
La facette $\mathcal{F}$ d\'efinie par $K$ est bien 
stable par $\theta$ car $K$ est d\'efini sur $F$. 

Prouvons que l'on a la condition b) de \ref{defbonneposition}.
Soit $R'$ le syst\`eme de racines de $G_L$,
relativement \`a $S_L$. Les \'el\'ements de $R$ et 
de $R'$ sont proportionnels (1.7. de \cite{[T79]}). 
Soit $R'_{\mathcal{F}}$ le syst\`eme des racines $\vec{\alpha}$
de $R'$ telles qu'il existe un mur de $\mathcal{A}_1$
de direction orthogonale \`a $\vec{\alpha}$ contenant 
$\mathcal{F}$. D'apr\`es 1.9. de \cite{[T79]},
$R'_{\mathcal{F}}$ s'identifie au syst\`eme de racines  
du quotient  $(\overline{K}_{\mathrm{red}})_k$ de la r\'eduction
modulo $\pi$ de $K_k$ par son radical unipotent
(\ref{g}), relativement au tore $\overline{T}_{\mathrm{red}}$
image de la r\'eduction de $S$ dans $\overline{K}_{\mathrm{red}}$.
Rappelons que $\overline{\mu}_{\mathrm{red}}$ est l'image 
de $\mu$ dans $\overline{K}_{\mathrm{red}}$ (prop. \ref{K}) ; 
c'est aussi l'image de $\mu_T$ (prop. \ref{calculinv}).
De plus, on a not\'e $\overline{T}_{\mu}$ le tore
de $\overline{K}_{\mathrm{red}}$ engendr\'e par 
les $\sigma ^i (\overline{\mu}_{\mathrm{red}})$,
$i\in \Z$. Notons $R'_{\mathcal{F},\lambda,\theta}$
le syst\`eme des racines de $R'_{\mathcal{F}}$ qui sont 
orthogonales aux cocaract\`eres de $(\overline{T}_{\mu})_k$.
Le syst\`eme de racines $R'_{\mathcal{F},\lambda,\theta}$
s'identifie \`a celui du commutant $(C_{\overline{T}_{\mu}})_k$ du tore
$(\overline{T}_{\mu})_k$ dans $(\overline{K}_{\mathrm{red}})_k$.
Comme le  tore maximal $\overline{T}_{\mathrm{red}}$ de
$C_{\overline{T}_{\mu}}$  contient un Borel de  
$C_{\overline{T}_{\mu}}$ (\ref{g}), 
l'action de $\sigma=\theta$
sur $R'_{\mathcal{F},\lambda,\theta}$ laisse stable une base
de $R'_{\mathcal{F},\lambda,\theta}$. Comme les \'el\'ements
de $R'_{\mathcal{F},\lambda,\theta}$ et de 
$R_{\mathcal{F},\lambda,\theta}$ sont proportionnels,
$\theta$ laisse stable une base de $R_{\mathcal{F},\lambda,\theta}$,
et on 
a bien la condition b) de \ref{defbonneposition}.

On peut donc 
choisir $\mathcal{C}_1$ telle que $(\mathcal{F},\lambda,\mathcal{C}_1)$  
soit en bonne position (prop. \ref{existebonneposition}).
On a donc prouv\'e que $X( \sigma(\{ \mu \}), b)_{\mathrm{Iw}_{\mathcal{C}}}$
est non vide. Il en r\'esulte que 
$X( \{ \mu \},\sigma^{-1}( b))_{\mathrm{Iw}_{\mathcal{C}}}$
est non vide. Comme $b$ et $\sigma^{-1}(b)$ sont $\sigma$-conjugu\'es,
on voit que $X( \{ \mu \}, b)_{\mathrm{Iw}_{\mathcal{C}}}$ est 
non vide. Comme $\mathrm{Iw}$ et $\mathrm{Iw}_{\mathcal{C}}$
sont conjugu\'es par un automorphisme int\'erieur de $G$,
on a bien que $X_w (b)_{\mathrm{Iw}}$ est non vide.
Ceci  prouve le th\'eor\`eme.

\subsection{}

Supposons $G$ non ramifi\'e. 
Si $K$ est un sous-groupe hypersp\'ecial de $G$,
$W^K$ s'identifie au groupe de Weyl $W_0$ de $G_{\overline{L}}$
et l'ensemble des doubles classes $W^K \backslash \widetilde{W} / W^K$
\`a $W_0 \backslash X_*$. En particulier, l'image de 
$\Lambda(\{ \mu \})$ dans $W^K \backslash \widetilde{W} / W^K$
est r\'eduit un \'el\'ement, que l'on note
encore $\{ \mu \}$. Si $\{ \mu \}$ est minuscule,
l'image de $ \mathrm{Adm}(\{ \mu\})$ dans $W_0 \backslash X_*$ est 
r\'eduite \`a $\{ \mu \}$ (prop. 3.11 de \cite{[R02]}).
On a donc le corollaire suivant, qui prouve la conjecture
4.6. de \cite{[R02]} lorsque $\{ \mu \}$ est minuscule.

\begin{cor}\label{conjhyper} Soit $K$ un sous-groupe hypersp\'ecial de $G$.
Soit $\{ \mu \}$ minuscule et $b\in B(G, \{ \mu \})$.
Alors $X_{ \{ \mu \}}(b)_K$ est non vide.\end{cor}

\subsection{Remarque.}\label{GP}
Soit $D$ un module de Dieudonn\'e filtr\'e
admissible et $M$ un r\'eseau fortement divisible. 
On suppose $\{ \mu \}$ minuscule \`a poids $0$ et $1$. 
Notons  $\overline{M}=M/pM$ 
avec sa structure de module de Dieudonn\'e filtr\'ee.
Cette structure d\'efinit sur $\overline{M}$  
une action du Frobenius $F$ et du Verschiebung $V=pF^{-1}$.  
Montrons 
comment on a proc\'ed\'e pour obtenir une structure iwahorique
\it i.e. \rm un drapeau de $\overline{M}$ stable par 
$F$ et $V$.
On a choisi tout d'abord une suite de  Jordan-H\"{o}lder
$(0)=\overline{M}_0\subset \overline{M}_1\subset
\ldots \overline{M}_r = \overline{M}$  du module de Dieudonn\'e filtr\'e
$\overline{M}$. Ensuite, pour chacun des quotients 
simples de cette suite de Jordan-H\"{o}lder, il n'est pas difficile
de voir que l' on a choisi dans le \ref{existebonneposition}
la filtration ``canonique'' du 2.5. de \cite{[Mo01]}
d\'efinie par $V$ et $F^{-1}$.

\nocite{*}
\bibliographystyle{plain}
\bibliography{kr5}

\def\cprime{$'$}
\begin{thebibliography}{10}

\bibitem{[SGA3]}
{\em Sch\'emas en groupes. {II}: {G}roupes de type multiplicatif, et structure
  des sch\'emas en groupes g\'en\'eraux}.
\newblock S\'eminaire de G\'eom\'etrie Alg\'ebrique du Bois Marie 1962/64 (SGA
  3). Dirig\'e par M. Demazure et A. Grothendieck. Lecture Notes in
  Mathematics, Vol. 152. Springer-Verlag, Berlin, 1962/1964.

\bibitem{[Bj83]}
Anders Bj{\"o}rner.
\newblock Orderings of {C}oxeter groups.
\newblock In {\em Combinatorics and algebra (Boulder, Colo., 1983)}, volume~34
  of {\em Contemp. Math.}, pages 175--195. Amer. Math. Soc., Providence, RI,
  1984.

\bibitem{[BT71]}
A.~Borel and J.~Tits.
\newblock \'{E}l\'ements unipotents et sous-groupes paraboliques de groupes
  r\'eductifs. {I}.
\newblock {\em Invent. Math.}, 12:95--104, 1971.

\bibitem{[LAG]}
Armand Borel.
\newblock {\em Linear algebraic groups}, volume 126 of {\em Graduate Texts in
  Mathematics}.
\newblock Springer-Verlag, New York, second edition, 1991.

\bibitem{[B98]}
Mikhail Borovoi.
\newblock Abelian {G}alois cohomology of reductive groups.
\newblock {\em Mem. Amer. Math. Soc.}, 132(626):viii+50, 1998.

\bibitem{[Bbki]}
N.~Bourbaki.
\newblock {\em \'{E}l\'ements de math\'ematique. {F}asc. {XXXIV}. {G}roupes et
  alg\`ebres de {L}ie. {C}hapitre {IV}: {G}roupes de {C}oxeter et syst\`emes de
  {T}its. {C}hapitre {V}: {G}roupes engendr\'es par des r\'eflexions.
  {C}hapitre {VI}: syst\`emes de racines}.
\newblock Actualit\'es Scientifiques et Industrielles, No. 1337. Hermann,
  Paris, 1968.

\bibitem{[BT72]}
F.~Bruhat and J.~Tits.
\newblock Groupes r\'eductifs sur un corps local.
\newblock {\em Inst. Hautes \'Etudes Sci. Publ. Math.}, (41):5--251, 1972.

\bibitem{[BT84]}
F.~Bruhat and J.~Tits.
\newblock Groupes r\'eductifs sur un corps local. {II}. {S}ch\'emas en groupes.
  {E}xistence d'une donn\'ee radicielle valu\'ee.
\newblock {\em Inst. Hautes \'Etudes Sci. Publ. Math.}, (60):197--376, 1984.

\bibitem{[FC01]}
Pierre Colmez and Jean-Marc Fontaine.
\newblock Construction des repr\'esentations {$p$}-adiques semi-stables.
\newblock {\em Invent. Math.}, 140(1):1--43, 2000.

\bibitem{[Fa95]}
Gerd Faltings.
\newblock Mumford-{S}tabilit\"at in der algebraischen {G}eometrie.
\newblock In {\em Proceedings of the International Congress of Mathematicians,
  Vol.\ 1, 2 (Z\"urich, 1994)}, pages 648--655, Basel, 1995. Birkh\"auser.

\bibitem{[FR02]}
J-M. Fontaine and R.~Rapoport.
\newblock Existence de filtrations admissibles sur des isocristaux.
\newblock {\em Preprint Mathematisches Institut des Universit\"{a}t zu
  K\"{o}ln}, pages 1--10, 2002.

\bibitem{[F94]}
Jean-Marc Fontaine.
\newblock Repr\'esentations {$p$}-adiques semi-stables.
\newblock {\em Ast\'erisque}, (223):113--184, 1994.
\newblock With an appendix by Pierre Colmez, P\'eriodes $p$-adiques
  (Bures-sur-Yvette, 1988).

\bibitem{[FL83]}
Jean-Marc Fontaine and Guy Laffaille.
\newblock Construction de repr\'esentations {$p$}-adiques.
\newblock {\em Ann. Sci. \'Ecole Norm. Sup. (4)}, 15(4):547--608 (1983), 1982.

\bibitem{[GP68]}
I.~M. Gel{\cprime}fand and V.~A. Ponomarev.
\newblock Indecomposable representations of the {L}orentz group.
\newblock {\em Russ. Math. Surv.}, 23:1--58, 1968.

\bibitem{[G98]}
Benedict~H. Gross.
\newblock Modular forms {$\pmod p$} and {G}alois representations.
\newblock {\em Internat. Math. Res. Notices}, (16):865--875, 1998.

\bibitem{[HR?]}
TH. Haines and R.~Rapoport.
\newblock Papier en pr\'eparation.

\bibitem{[IM65]}
N.~Iwahori and H.~Matsumoto.
\newblock On some {B}ruhat decomposition and the structure of the {H}ecke rings
  of {${ p}$}-adic {C}hevalley groups.
\newblock {\em Inst. Hautes \'Etudes Sci. Publ. Math.}, (25):5--48, 1965.

\bibitem{[KR02]}
R.~Kottwitz and M.~Rapoport.
\newblock On the existence of {${F}$}-crystals.
\newblock {\em arXiv:math.NT/0202229}, pages 1--28, 2002.

\bibitem{[K85]}
Robert~E. Kottwitz.
\newblock Isocrystals with additional structure.
\newblock {\em Compositio Math.}, 56(2):201--220, 1985.

\bibitem{[K97]}
Robert~E. Kottwitz.
\newblock Isocrystals with additional structure. {I}{I}.
\newblock {\em Compositio Math.}, 109(3):255--339, 1997.

\bibitem{[Kraft]}
H~Kraft.
\newblock Kommutative algebraische $p$-{G}ruppen (mit {A}nwendungen auf
  $p$-divisible {G}ruppen und abelsche {V}ariet\"aten).
\newblock {\em Manuscrit non publi\'e}, 1975.

\bibitem{[L80]}
Guy Laffaille.
\newblock Groupes {$p$}-divisibles et modules filtr\'es: le cas peu ramifi\'e.
\newblock {\em Bull. Soc. Math. France}, 108(2):187--206, 1980.

\bibitem{[Le03]}
C.~Leigh.
\newblock A converse to {M}azur's inequality for split clasical groups.
\newblock {\em arXiv:math.NT/0211327}, pages 1--16, 2002.

\bibitem{[M72]}
B.~Mazur.
\newblock Frobenius and the {H}odge filtration.
\newblock {\em Bull. Amer. Math. Soc.}, 78:653--667, 1972.

\bibitem{[Mo01]}
Ben Moonen.
\newblock Group schemes with additional structures and {W}eyl group cosets.
\newblock In {\em Moduli of abelian varieties (Texel Island, 1999)}, volume 195
  of {\em Progr. Math.}, pages 255--298. Birkh\"auser, Basel, 2001.

\bibitem{[R02]}
M~Rapoport.
\newblock A guide to the reduction modulo $p$ of {S}himura varieties.
\newblock {\em Preprint}, 02.

\bibitem{[RR96]}
M.~Rapoport and M.~Richartz.
\newblock On the classification and specialization of {$F$}-isocrystals with
  additional structure.
\newblock {\em Compositio Math.}, 103(2):153--181, 1996.

\bibitem{[RZ96]}
M.~Rapoport and Th. Zink.
\newblock {\em Period spaces for {$p$}-divisible groups}, volume 141 of {\em
  Annals of Mathematics Studies}.
\newblock Princeton University Press, Princeton, NJ, 1996.

\bibitem{[R01]}
Michael Rapoport.
\newblock A positivity property of the {S}atake isomorphism.
\newblock {\em Manuscripta Math.}, 101(2):153--166, 2000.

\bibitem{[Re03]}
D.~C. Reuman.
\newblock Determining whether certain affine {D}eligne-{L}usztig sets are
  empty.
\newblock {\em arXiv:math.NT/0211434}, pages 1--135, 2002.

\bibitem{[Sa72]}
Neantro Saavedra~Rivano.
\newblock {\em Cat\'egories {T}annakiennes}.
\newblock Springer-Verlag, Berlin, 1972.
\newblock Lecture Notes in Mathematics, Vol. 265.

\bibitem{[S79]}
Jean-Pierre Serre.
\newblock Groupes alg\'ebriques associ\'es aux modules de {H}odge-{T}ate.
\newblock In {\em Journ\'ees de G\'eom\'etrie Alg\'ebrique de Rennes. (Rennes,
  1978), Vol. III}, volume~65 of {\em Ast\'erisque}, pages 155--188. Soc. Math.
  France, Paris, 1979.

\bibitem{[CG]}
Jean-Pierre Serre.
\newblock {\em Cohomologie galoisienne}, volume~5 of {\em Lecture Notes in
  Mathematics}.
\newblock Springer-Verlag, Berlin, fifth edition, 1994.

\bibitem{[T79]}
J.~Tits.
\newblock Reductive groups over local fields.
\newblock In {\em Automorphic forms, representations and $L$-functions (Proc.
  Sympos. Pure Math., Oregon State Univ., Corvallis, Ore., 1977), Part 1},
  pages 29--69. Amer. Math. Soc., Providence, R.I., 1979.

\bibitem{[W84]}
Jean-Pierre Wintenberger.
\newblock Un scindage de la filtration de {H}odge pour certaines vari\'et\'es
  alg\'ebriques sur les corps locaux.
\newblock {\em Ann. of Math. (2)}, 119(3), 1984.

\bibitem{[W97]}
Jean-Pierre Wintenberger.
\newblock Propri\'et\'es du groupe tannakien des structures de {H}odge
  {$p$}-adiques et torseur entre cohomologies cristalline et \'etale.
\newblock {\em Ann. Inst. Fourier (Grenoble)}, 47, 1997.

\end{thebibliography}

\vspace{2cm}
Jean-Pierre Wintenberger

Universit\'e Louis Pasteur

D\'epartement de Math\'ematiques, IRMA

7, rue Ren\'e Descartes 

67084 Strasbourg Cedex

France

e-mail wintenb@math.u-strasbg.fr

tel 03 90 24 02 17 , fax 03 90 24 03 28

\end{document}